\documentclass{amsproc}

\newtheorem{theorem}{Theorem}[section]
\newtheorem{lemma}[theorem]{Lemma}

\theoremstyle{definition}
\newtheorem{definition}[theorem]{Definition}

\theoremstyle{remark}

\numberwithin{equation}{section}

\usepackage{amssymb}
\usepackage{amsmath}
\usepackage{latexsym}
\usepackage[mathscr]{euscript}
\usepackage{graphicx}
\usepackage{amscd}
\usepackage{amsthm} 
\usepackage{fullpage} 
\usepackage{wrapfig}
\newtheorem{corollary}[theorem]{Corollary}
\newtheorem{proposition}[theorem]{Proposition}

\newtheorem{induct}[theorem]{Inductive Assumption}

\newcommand{\bc}{\mathbb{C}}
\newcommand{\bp}{\mathbb{ P}}
\newcommand{\bz}{\mathbb{Z}}
\newcommand{\br}{\mathbb{R}}
\newcommand{\bq}{\mathbb{Q}}

\newcommand{\p}{\partial}
\newcommand{\cc}{\mathcal{C}}

\newcommand{\hk}{\hookrightarrow}
\newcommand{\bg}{\bigskip}
\newcommand{\med}{\medskip}
\newcommand{\la}{\longrightarrow}
\newcommand{\bfl}{\begin{flushleft}}
\newcommand{\efl}{\end{flushleft}}

\newcommand{\eps}{\epsilon}
\newcommand{\hocolim}{\operatorname{hocolim}}

\newcommand{\xr}{\xrightarrow}

\newcommand{\dgn}{Diff(F_{g,n}, \p)}

\newcommand{\si}{\Sigma_{1,0}}
\newcommand{\sj}{\Sigma_{0,1}}
\newcommand{\sk}{\Sigma_{0,-1}}
\newcommand{\sij}{\Sigma_{i,j}}
\newcommand{\di}{\Delta_{1,0}}
\newcommand{\djj}{\Delta_{0,1}}

\newcommand{\dij}{\Delta_{i,j}}
\newcommand{\sgn}{\mathscr S_{g,n}}
\newcommand{\fgn}{F_{g,n}}
\newcommand{\sinf}{\mathscr S_{\infty, n}}
\newcommand{\G}{\Gamma}
\newcommand{\msc}{\mathscr}
\newcommand{\rv}{Rel^V}
\newcommand{\hf}{\hat F}

 \newcommand{\fgp}{F_{g-p, 2p+1}}
 \newcommand{\Kappa}{\msc K}
 \newcommand{\cred}{\cc^{red}_X}
 \newcommand{\sst}{(S^1, *)}
 \newcommand{\tst}{(T, *)}

\begin{document}  

  \title{Surfaces in a background space and the homology of mapping class groups}
  \author{Ralph L. Cohen} 
  \address{Department of Mathematics, Bldg. 380, Stanford University, Stanford, CA 94305, USA} 
  \email{ralph@math.stanford.edu}
   \thanks{Both authors were partially supported by a Focused Research Group grant  from the NSF}
\author {Ib Madsen}
 \address{ Department of Mathematical Sciences,  University of  Copenhagen,  Universitetsparken 5, DK-2100 Copenhagen, Denmark }
 \email{imadsen@imf.au.dk}
\date{\today}

\maketitle  
 \begin{abstract}
In this paper we study the topology of the space of Riemann surfaces in a simply connected space $X$,  $\sgn (X, \gamma)$.  This is the space consisting of triples, $(F_{g,n}, \phi, f)$, where $F_{g,n}$ is a Riemann surface of genus $g$ and $n$-boundary components, $\phi$ is a parameterization of the boundary, $\p F_{g,n}$, and $f : F_{g,n} \to X$ is a continuous map that satisfies a boundary condition $\gamma$.  We prove three theorems about these spaces.  Our main theorem is the identification of the stable homology type of the space $\msc S_{\infty, n}(X; \gamma)$, defined to be the limit  as the genus $g$ gets large, of the spaces  $\sgn (X; \gamma)$.  Our result about this stable topology is a parameterized version of the theorem of Madsen  and Weiss proving a generalization of the Mumford conjecture on the stable cohomology of mapping class groups.  Our second result describes a stable range in which the homology of $\sgn (X; \gamma)$ is isomorphic to the stable homology.  Finally we prove a stability theorem about the homology of mapping class groups with certain families of twisted coefficients.  The second and third theorems are generalizations of stability theorems of Harer and Ivanov.   \end{abstract}

 \tableofcontents

 \section*{Introduction}

 The goal of this paper is to study the topology of the space of surfaces mapping to a background space $X$, with boundary condition $\gamma$,   $\sgn (X; \gamma)$.
 This space is defined as follows.  
 
Let $X$ be a   simply connected space with basepoint $x_0 \in X$.   Let $\gamma : \coprod_n S^1 \to X$ be $n$ continuous loops in $X$.  Define the space
\begin{align}
\sgn (X, \gamma) = &\{(S_{g,n}, \phi, f): \,  \text{where $S_{g,n} \subset \br^\infty \times [a,b]$ is a smooth oriented surface of genus $g$ and } \notag  \\
 &\text{ $n$ boundary components, $\phi : \coprod_n S^1 \xr{\cong} \p S$ is a parameterization of the boundary,}\notag \\
& \text{   and $f : S_{g,n} \to X$ is a continuous map with $\p f = \gamma : \coprod_n S^1 \to X$. } \} \notag
\end{align}
In this description, $[a,b]$ is an arbitrary closed interval, and the boundary, $\p S$, lies in the boundary, $\p S = (\br^\infty \times \{a\}  \sqcup \br^\infty \times \{b\})$. We also insist that if $n>0$, the ``incoming" boundary, $\p S \cap (\br^\infty \times \{a\})$ has one connected component, which we refer to as $\p_0 S$.  The parameterization $\phi$ is an orientation preserving diffeomorphism.        $\p f$ is the composition $\coprod_n S^1 \xr{\phi} \p S \xr{f_{|_{\p S}}} X$. 

\med
We think of these spaces as moduli spaces of  Riemann surfaces mapping to $X$, or  for short, the moduli space of surfaces in $X$.   Indeed the embedding of the surface in Euclidean space defines an inner product on the tangent space of the surface, which together with the orientation defines an almost complex structure, and hence a complex structure on the surface. 

\med
We have three main results in this paper.  The first describes the ``stable topology" of $\sgn (X, \gamma)$, the second is a stability result showing the range of dimensions in which the homology of $\sgn (X; \gamma)$ is in the stable range, and the third is a stability  result about the homology of mapping class groups with certain families of twisted coefficients.  

We haven't yet described the topology of $\sgn (X; \gamma)$.  To do this, let $F_{g,n}$ be a fixed surface of genus $g$ with $n$ boundary components.  Let $\delta : \coprod_n S^1 \xr{\cong} \p \fgn$  be a fixed parameterization of the boundary.  Let $Emb(\fgn, \br^\infty)$ be the space of embeddings $e : F_{g,n} \hk  \br^\infty \times [a,b]$ as above, for some choice of $a < b$.   The topology on this space is induced by the compact open topology.  The Whitney embedding theorem implies  that  $Emb(\fgn, \br^\infty)$ is   contractible.     It also has    a free action of the group $\dgn$ of orientation preserving  diffeomorphisms
of $\fgn$ that fix the boundary pointwise.  The action is given by precomposition.  Let $Map_\gamma (\fgn, X)$ be the space of continuous maps $f : \fgn \to X$ with $\p f = \gamma$.  This also has the compact-open topology.  It is also acted up by $\dgn$ by precomposition.   We then have the following immediate observation.

\bfl
\bf
Observation. \rm There is a bijective correspondence,
\begin{align}
\sgn (X; \gamma) &\cong  Emb(\fgn, \br^\infty) \times_{\dgn}Map_\gamma (F_{g,n}, X) \notag \\
&\simeq E(\dgn) \times_{\dgn} Map_\gamma (F_{g,n}, X) \notag 
\end{align}

\efl

Notice in particular that when $X$ is a point,  the space of surfaces is the classifying space of the diffeomorphism group, $\sgn (point) \simeq BDiff(\fgn).$

\med
We next observe that the spaces $Map_\gamma (F_{g,n}, X)$ and $\sgn (X; \gamma)$ have homotopy types that do not depend on the boundary map $\gamma$. This is for the following reason.  Consider the mapping spaces, $Map (F_{g,n}, X)$ and $\sgn (X)$ that have no boundary conditions at all.  Then restriction of these mapping spaces to the boundary, determines Serre fibrations,
$$
Map (F_{g,n}, X) \to (LX)^n    \quad \text{and} \quad \sgn (X) \to (LX)^n
$$
where $LX = Map (S^1, X)$ is the free loop space.  Since $X$ is assumed to be simply connected the base spaces of these fibrations, $(LX)^n$, are connected.  Therefore the fibers of these maps have homotopy types which are independent of the choice of point
$\gamma \in (LX)^n$.  

Because of this fact, we are free to work with convenient choices of boundary conditions.
We will assume our boundary map $\gamma : \coprod_n S^1 \to X$, viewed as $n$-loops numbered $\gamma_0, \cdots , \gamma_{n-1}$, has the  property that $\gamma_0 : S^1 \to x_0 \in X$ is constant at the basepoint.

Notice that given a point $(S, \phi, f) \in \sgn (X; \gamma)$, the above numbering
and the parameterization $\phi$ determines a numbering the boundary components, $\p_0S, \cdots , \p_{n-1} S$.  
Also  the boundary components of $S$ are partitioned as a disjoint union, $\p S = \p_a S \sqcup \p_b S$, the ``incoming" and ``outgoing" components of the boundary. We assume that $\p_0 S \in \p_a S$ is an incoming boundary component.    By the boundary conditions, $\p_0 S$ is mapped by $f$  to the basepoint  $x_0\in X$.

   The boundary components are oriented in two different ways, namely by the parameterization $\phi$, and by the induced orientation from $S$.  We assume that the two orientations are opposite for the incoming components, $\p_aS$, and agree for the outgoing components $\p_bS$.   
   
\med
To state our result about the stable topology  of $\sgn (X; \gamma)$,  fix a  surface of genus one,   $T \subset \br^3 \times [0,1] \subset \br^\infty \times [0,1]$ having one incoming and one outgoing boundary component.

\begin{figure}[ht]
  \centering
  \includegraphics[height=4cm]{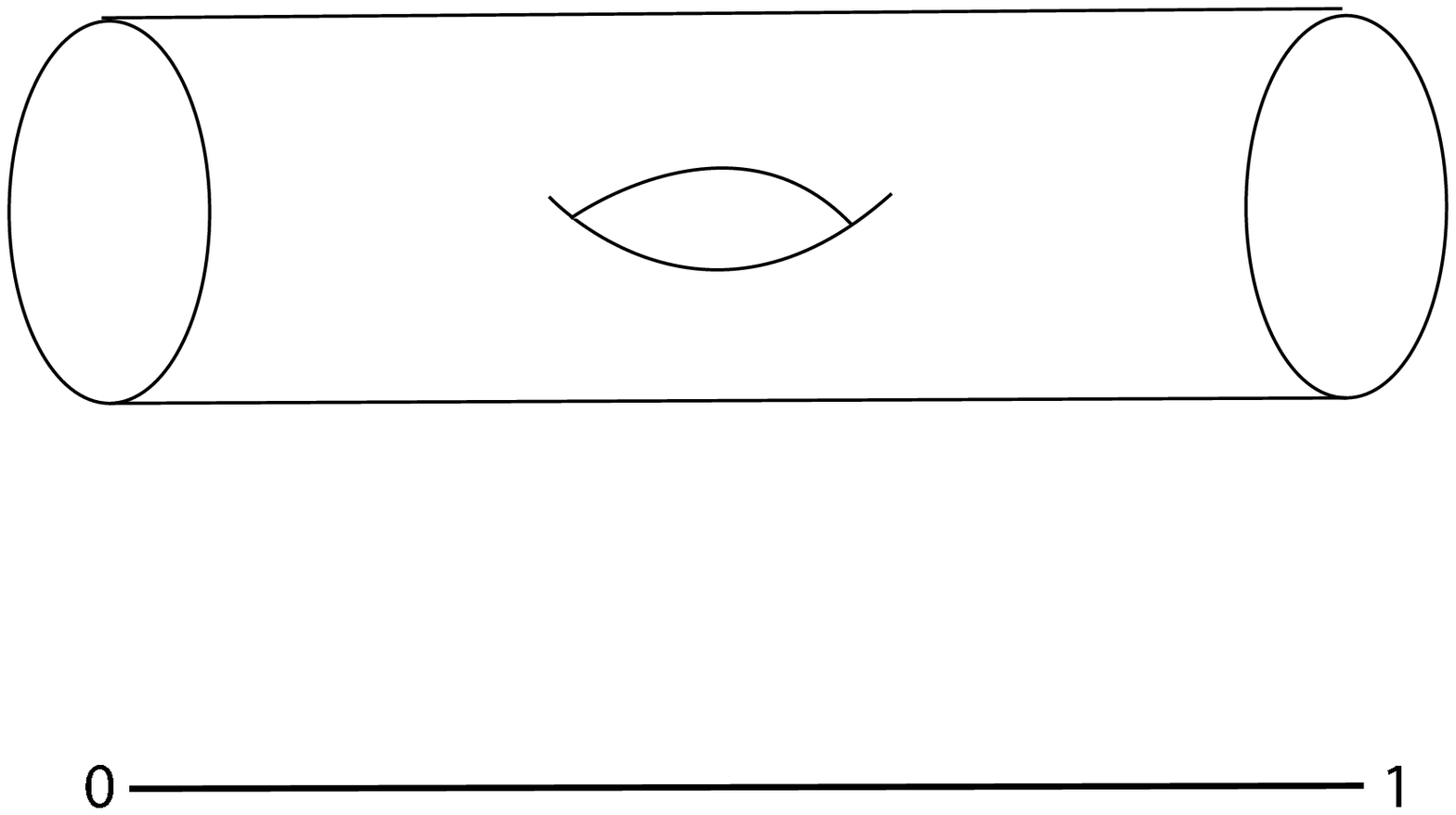}
  \caption{The surface $T \subset \br^\infty \times [0,1]$}
  \label{fig:figone}
\end{figure}

Given $(S, \phi, f) \in \sgn (X, \gamma)$, we  ``glue in" the surface $T$ to get an element of $\mathscr S_{g+1,n}(X; \gamma)$ as follows.  Suppose $S \subset \br^\infty \times [a,b]$.  Translate $T$ so it is now embedded in $\br^\infty \times [a-1, a]$.  Identify the boundary $\p_a T$ with $\p_0 S$ using the parameterizations.  Similarly glue in a cylinder $S^1 \times [a-1,a]$ to each of the other boundary components in $\p_a S$.  The result is a surface $T\#S$ of genus $g+1$ embedded in $\br^\infty \times [a-1, b]$.  
The boundary  parameterization $\phi$ now defines a boundary parameterization  of $T\# S$, and the map $f : S \to X$ extends to $T\# S$ by letting it be constant at the basepoint on $T$, and on each new  cylinder glued in on the $i^{th}$ boundary $\p_i S \subset \p_a S$, it is defined to be the composition $S^1 \times [a-1,a] \xr{project} S^1 \xr{\p_i f} X$.   This construction defines a map
$$
T_\# : \sgn (X; \gamma) \to \mathscr S_{g+1,n} (X; \gamma).
$$

\begin{figure}[ht]
  \centering
  \includegraphics[height=6cm]{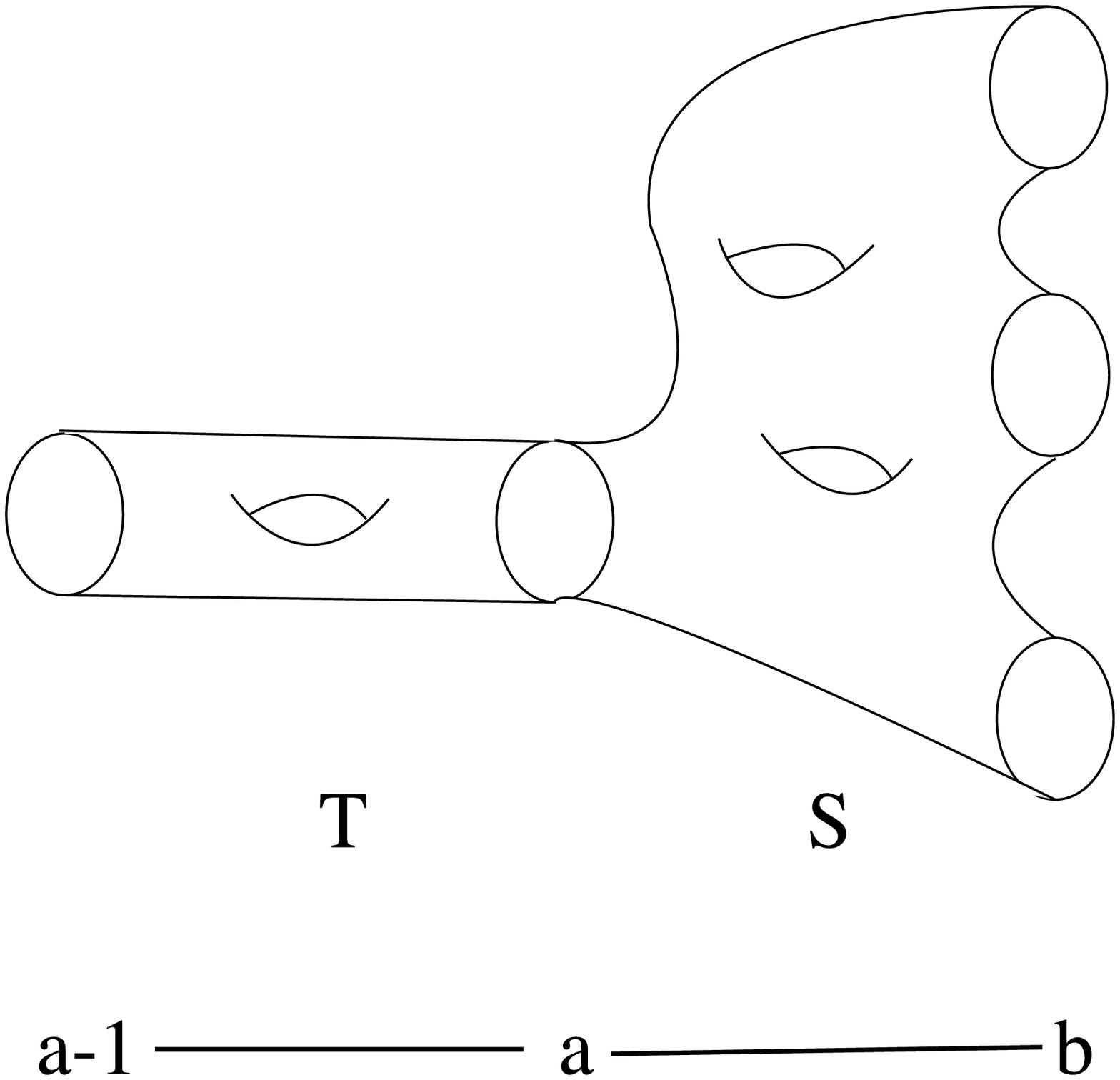}
  \caption{  $T_ \#S$}
  \label{fig:figtwo}
\end{figure}

We now define $\sinf (X; \gamma)$ to be the homotopy colimit of the map $T_\#$,
$$
\sinf (X; \gamma) =  hocolim  \,\{\sgn (X; \gamma) \xr{T_\#} \mathscr S_{g+1,n} (X; \gamma)  \xr{T_\#} \cdots \}
$$
Recall that in this situation the homotopy colimit is the infinite mapping cylinder
of the iterations of the map $T_\#$.  Thus it is a particular type of direct, or colimit.

We refer to the topology of $\sinf (X; \gamma)$ as the ``stable topology" of the moduli spaces,  $\sgn (X; \gamma)$.

 Our first theorem describes the stable topology of these moduli spaces.
 
 \med
 \begin{theorem}\label{one} Let $X$ be a simply connected, based space.  There is a map,
 $$
\alpha :  \bz \times \sinf (X; \gamma) \la \Omega^\infty (\bc \bp^\infty_{-1} \wedge X_+).
 $$ that induces an isomorphism in any generalized homology theory. 
 \end{theorem}
 In this theorem, the right hand side is the infinite loop space defined to be the zero space of the spectrum $\bc \bp^\infty_{-1} \wedge X_+$.  Here $\bc \bp^\infty_{-1}$
 is the Thom spectrum of the virtual bundle $-L \to \bc \bp^\infty$, where $L \to \bc \bp^\infty$ is the canonical line bundle over $\bc \bp^\infty$. 
 
 \med
 We observe that when $X$ is a point,  $\sinf (point)$ represents, up to homotopy, the stable topology of the moduli
 space of bordered Riemann surfaces studied by the second author and Weiss in \cite{madsenweiss}, in their proof of the generalized Mumford conjecture.  Indeed, the Madsen-Weiss theorem is a key ingredient in our proof of Theorem \ref{one}, and in the case when $X = point$,   Theorem  \ref{one} is just a restatement of their theorem.  Thus Theorem \ref{one}
 can be viewed as a parameterized form of the Madsen-Weiss theorem, where $X$ is the parameterizing space. 
 
 \med
 We remark that the homology of the infinite loop space in this theorem has been completely computed by Galatius \cite{galatius} when $X$ is a point.  The rational cohomology is much simpler.    The following corollary states that the rational stable cohomology of the space of surfaces in $X$ is     generated by the Miller-Morita-Mumford  $\kappa$-classes, and the rational cohomology of $X$.
 
 To state this more carefully, we restrict our attention to a particular path component of $ \sinf (X; \gamma)$.  It is clear that since $X$ is simply connected,  the set of path components $\pi_0( \sinf (X; \gamma))$ is in bijective correspondence with the path components
 of the mapping space, $\pi_0(Map_\gamma (F_{g,n}, X))$, which in turn is in bijective correspondence with the homotopy group,
 $\pi_2(X)$.  Moreover, since Theorem \ref{one} tells us that  $ \sinf (X; \gamma)$ is homology equivalent to an infinite loop space, all of its path components have isomorphic homologies.  Let $ \sinf (X; \gamma)_\bullet$ be the connected path component corresponding to the trivial class $0 \in \pi_2(X)$.  Similarly let $\Omega^\infty_\bullet (\bc \bp^\infty_{-1} \wedge X_+)$ represent the corresponding
connected path component.   Since $H^*(\sinf (X; \gamma)_\bullet ; \bq) \cong H^*(\Omega^\infty_\bullet (\bc \bp^\infty_{-1} \wedge X_+) ; \bq)$, then \cite{milnormoore} gives us the following description of the rational cohomology.

 Suppose $V$ is a graded vector space over the rationals, and $A(V)$ is the free $\bq$-algebra generated by $V$.   That is, given a basis of $V$, $A(V)$ is the polynomial algebra generated by the even dimensional basis elements, tensor the exterior algebra generated by the odd dimensional basis elements. 
 
 Let $\Kappa$ be the graded vector space over $\bq$ generated by one basis element, $\kappa_i$, of dimension $2i$ for each $i \geq -1$.  Consider the tensor product of graded vector spaces, $\Kappa \otimes H^*(X; \bq)$.  Let $ (\Kappa \otimes H^*(X; \bq))_+$ denote
 that part of this vector space that lives in positive grading.  We then have the following.

 \med
 \begin{corollary}  There is an isomorphism of algebras,
 $$
 H^*(\sinf (X; \gamma)_\bullet \, \, \bq) \cong A((\Kappa \otimes H^*(X; \bq))_+). $$
 \end{corollary}
  
 As we remarked before, $H^*(\sinf (point); \bq)$ is the stable rational cohomology of moduli space.  This algebra was conjectured by Mumford, and proven by Madsen and Weiss in     \cite{madsenweiss}, to be    the polynomial algebra on the Miller-Morita-Mumford $\kappa$-classes.  The classes $\kappa_i \in \Kappa \subset  H^*(\sinf (X; \gamma) \, \bq)$ for $i \geq 1$ are the image of the Miller-Morita-Mumford  classes under the map $$H^*(\sinf (point); \bq) \to H^*(\sinf (X; \gamma) \, \bq).$$

 \med
 Now notice that in the statement of Theorem \ref{one}, the right hand side does not depend on $n$, the number of boundary components.     This is strengthened by the following theorem,
 which identifies the stable range of the homology of the individual surface spaces.
 
 \med
 \begin{theorem}\label{unstable}  For $X$ simply connected as above,  the homology groups,
 $$
 H_q(\sgn(X; \gamma))
 $$
 are independent of the genus $g$, the number of boundary components $n$, and the boundary condition $\gamma$, so long as $2q+4 \leq g$.  In other words, for $q$ in this range,
 $$H_q(\sgn(X; \gamma)_\bullet) \cong H_q(\sinf(X; \gamma)_\bullet) \cong H_q(\Omega^\infty_\bullet (\bc \bp^\infty_{-1} \wedge X_+)).
 $$
 \end{theorem}

 \med
 Our last result,  which is actually a key ingredient in proving both Theorem \ref{one} and  Theorem \ref{unstable} \, is purely a statement about the homology of groups.  Our inspiration for this theorem was the work of Ivanov \cite{ivanov} which gave the first stability results
 for the homology of mapping class groups with certain kinds of twisted coefficients.
 The  following is a generalization of his results.   
 
 Let $\G_{g,n} = \pi_0(\dgn)$ be the mapping class group. Notice there are natural maps,
 $$
T\# : \G_{g,n} \to  \G_{g+1,n}     \quad \text{and} \quad P\# :\G_{g,n} \to  \G_{g ,n+1} 
   $$
   induced by gluing in the surface of genus one, $T$, as above, and by gluing in a ``pair of pants" $P$, of genus zero, with two incoming and one outgoing boundary component.
   The gluing procedure is completely analogous to the gluing in of the surface $T$ described above.

\begin{figure}[ht]
  \centering
  \includegraphics[height=4cm]{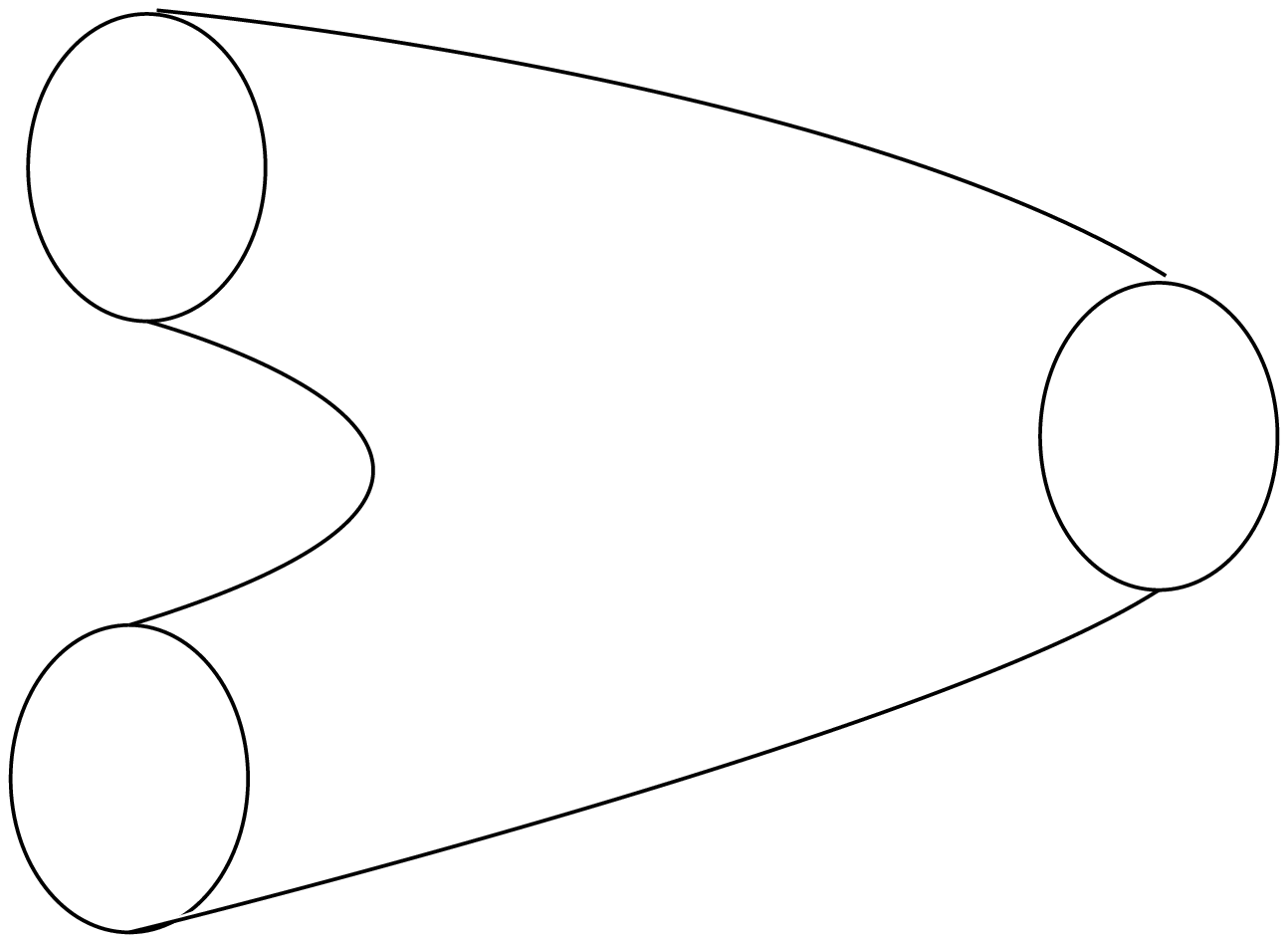}
  \caption{The pair of pants surface $P$}
  \label{fig:figthree}
\end{figure}

  In section 1 we will define the notion of a    \sl coefficient system \rm $V = \{V_{g,n}\}$,  generalizing the notion defined by Ivanov \cite{ivanov}. However for the purposes of the statement of the following theorem, its main property is    the following.  $V$  is a collection of modules
   $V_{g,n}$ over $\bz[\G_{g,n}]$, together with split injective maps     $$
   \si : V_{g,n} \to V_{g+1,n } \quad \text{and} \quad \sj :  V_{g,n} \to V_{g,n+1}$$ that are equivariant with respect to the gluing maps $T\#$ and $P\#$ respectively, so that we have splittings
   \begin{align}
   V_{g+1,n} &\cong V_{g,n} \oplus \Delta_{1,0}V_{g,n} \notag \\
   V_{g, n+1} &\cong V_{g,n} \oplus \Delta_{0,1}V_{g,n} \notag
   \end{align}
   as $\G_{g,n}$-modules.   We say that the coefficient system $V$ has \sl degree \rm $0$ if it is constant; that is all the modules $\Delta_{1,0}V_{g,n}$ and $\Delta_{0,1}V_{g,n}$ are zero.  Recursively, we define the degree of $V$ to be $d$, if the coefficient systems $\{\Delta_{1,0}V_{g,n}\}$ and $\{\Delta_{0,1}V_{g,n}\}$ have degree $\leq d-1$.  
   
   A nice example of a coefficent system of degree  one  is $V_{g,n}= H_1(F_{g,n} ; \bz)$.    Our last main result  is a generalization of stability theorems of Harer \cite{harer} and Ivanov \cite{ivanov}, \cite{ivanov2}.
   
   \med
   \begin{theorem}\label{mpg}
   If $V$ is a coefficient system of degree $d$, then the   homology group
   $$
   H_q(\G_{g,n}; V_{g,n})
   $$ 
   is independent of $g$, and $n$, so long as $2q+d +2< g-1$.  If we require $n$ to be positive, then this stability range improves to $2q+d +2< g$.
   \end{theorem}

   \med
   This paper is organized as follows.  In section one we will prove Theorem \ref{mpg}.
   This will involve adaptation and generalization of ideas of Ivanov  \cite{ivanov}.  
   In section 2 we use this result to prove Theorem \ref{unstable}.  To do this we use Theorem \ref{mpg} and   a homotopy theory argument and calculation.  Finally in section 3 we prove Theorem \ref{one}.  We do this by adapting     arguments of Tillmann \cite{tillmann} using techniques of McDuff-Segal \cite{mcduffsegal},  to show how Theorem \ref{one} ultimately follows from the Madsen-Weiss theorem \cite{madsenweiss} and Theorem \ref{unstable}. 
   
   \med
   Both authors owe a debt of gratitude to Ulrike Tillmann for many hours of very helpful conversation
   about the arguments and results of this paper.  We also thank Nathalie Wahl for help with the curve complexes used below, as well as Elizabeth Hanbury who spotted an error in the definition of our surface category in a previous version.

   \section{The homology of mapping class groups with twisted coefficients}
   Our goal in this section  is to prove Theorem \ref{mpg}, as stated in the introduction.  To do this we use ideas of Ivanov \cite{ivanov},  generalized to the context that we need.

   
  \subsection{Categories of surfaces and coefficient systems}
   
   We begin by defining the category of differentiable surfaces $\cc$ in which we will work.
  
  \begin{definition}\label{category}
  For $g, n \geq 0$, we define the category $\cc_{g,n}$ to have   objects $(F, \phi)$,  where $F$ is a smooth, oriented, compact surface,  of genus $g$ and $n$  boundary components,  and $\phi : \coprod_n S^1 \xr{\cong} \p F$ is an orientation preserving diffeomorphism (i.e a  parameterization)  of the boundary.  
  We write $\phi = \phi_0, \cdots , \phi_{n-1}$,  which has the effect of numbering the boundary components of $F$, $\p_0F, \cdots \p_{n-1}F$. 
 
 A  morphism    $e : (F_1, \phi_1) \to (F_2, \phi_2)$ is an   isotopy class of orientation preserving  diffeomorphism $F_1 \to  F_2$, that preserves the boundary parameterizations.       \end{definition}
 
     \med
   We now put all these categories together.   If $(F, \phi)$ is an object in $\cc_{g,n}$ with $n \geq 1$, let   $x_0 \in \p_0F$ be the basepoint   that corresponds to the basepoint $0 \in \br/\bz = S^1$ under the parameterization $\phi$. 
   
      \begin{definition}\label{fullcat}   
     Define the surface category $\cc$  to have   objects equal to  the disjoint union,
     $$
     Ob \, (\cc) = \coprod_{g,n} Ob \,  \cc_{g,n}.$$
     There is a morphism for each ambient isotopy class of embedding,   $e : F_1 \to F_2$,  that maps each boundary component of $F_1$ either diffeomorphically to a boundary component of $F_2$,  respecting the parameterizations, or to the interior of $F_2$.   If the boundary $\p F_2 \neq \emptyset$, and if $e : F_1 \hk F_2$ maps the boundary component $\p_0 F_1$ to the interior
     of $F_2$,    then we also require in our definition of morphism 
  an ambient isotopy class of  parameterized embedded arc $\gamma : [0,t] \hk  F_2$, for some $t \geq 0$,   starting at the basepoint    $e(x_0) \in e(\p_0 F_1)$ and ending at the basepoint $x_0 \in \p_0 F_2$.   The interior of the embedding $\gamma$ is required to lie in the interior of $F_2$. 
     \end{definition}

\begin{figure}[ht]
  \centering
  \includegraphics[height=6cm]{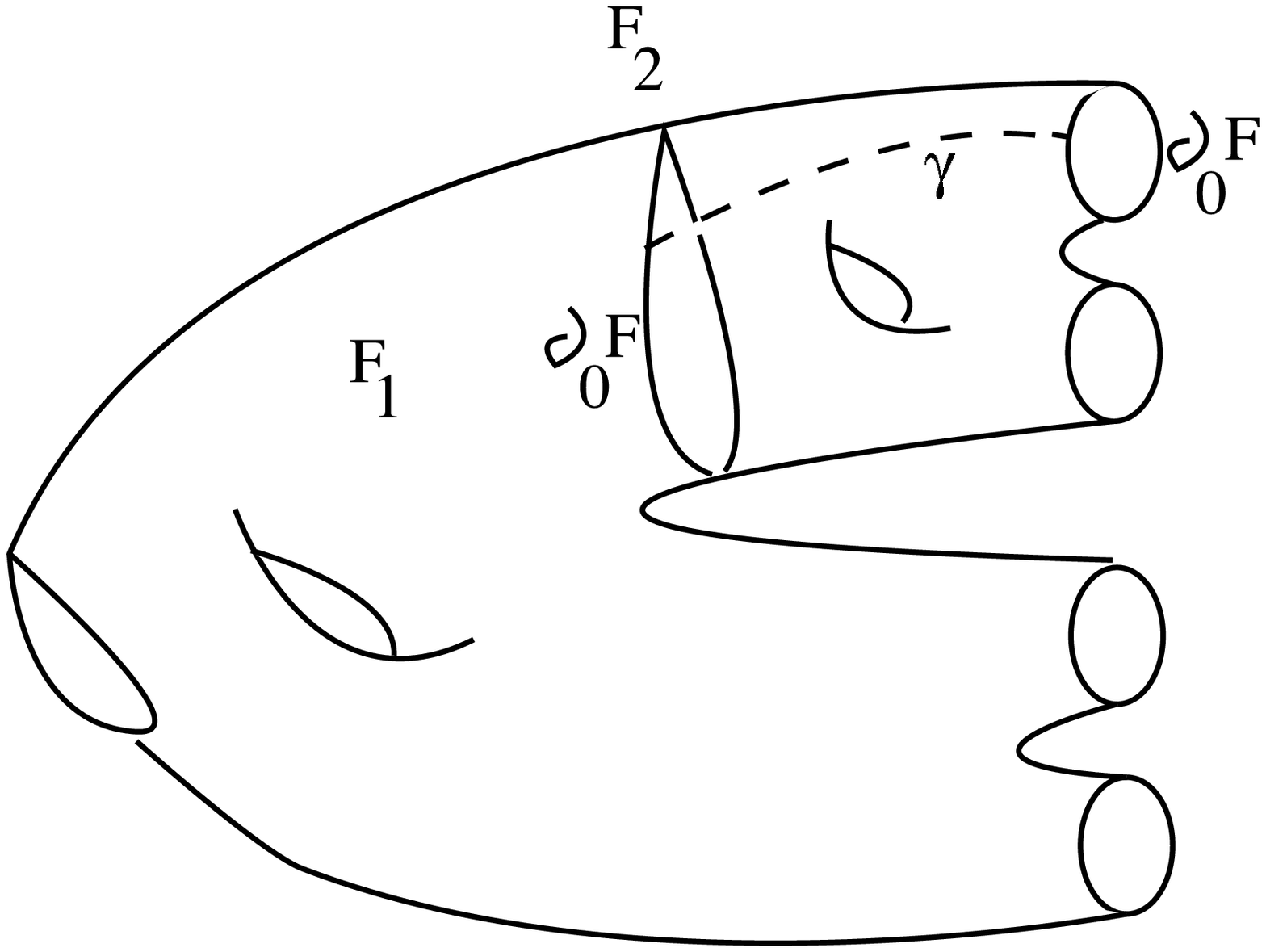}
 \caption{ A morphism $e : F_1 \to  F_2$  } 
  \label{morphism}
\end{figure}

     We remark that this category is a slight variation of     Ivanov's  category  of decorated surfaces in \cite{ivanov}.   Notice that each $\cc_{g,n} \subset \cc$ is a subcategory.

     \med
       We now describe three    functors  $\si : \cc_{g,n} \to \cc_{g+1,n}$,     $\sj : \cc_{g,n} \to \cc_{g,n+1}$, and $\sk: \cc_{g,n} \to \cc_{g ,n-1}$  defined when $n \geq 1$.      These operations  have the 
   effect of increasing the genus by one,   increasing the number of boundary components by one, and decreasing the number of boundary components by one, respectively. 
   
   In order to describe these operations more precisely, we use a graphical technique of Ivanov \cite{ivanov}.  Namely, since these operations are all defined on surfaces with at least one boundary component, then we can  identify $\p_0F$  with a rectangle using the parameterization, and we can picture  $F$  as a rectangular disk with handles attached, and disks removed  from the interior.

   \begin{figure}[ht]
  \centering
  \includegraphics[height=6cm]{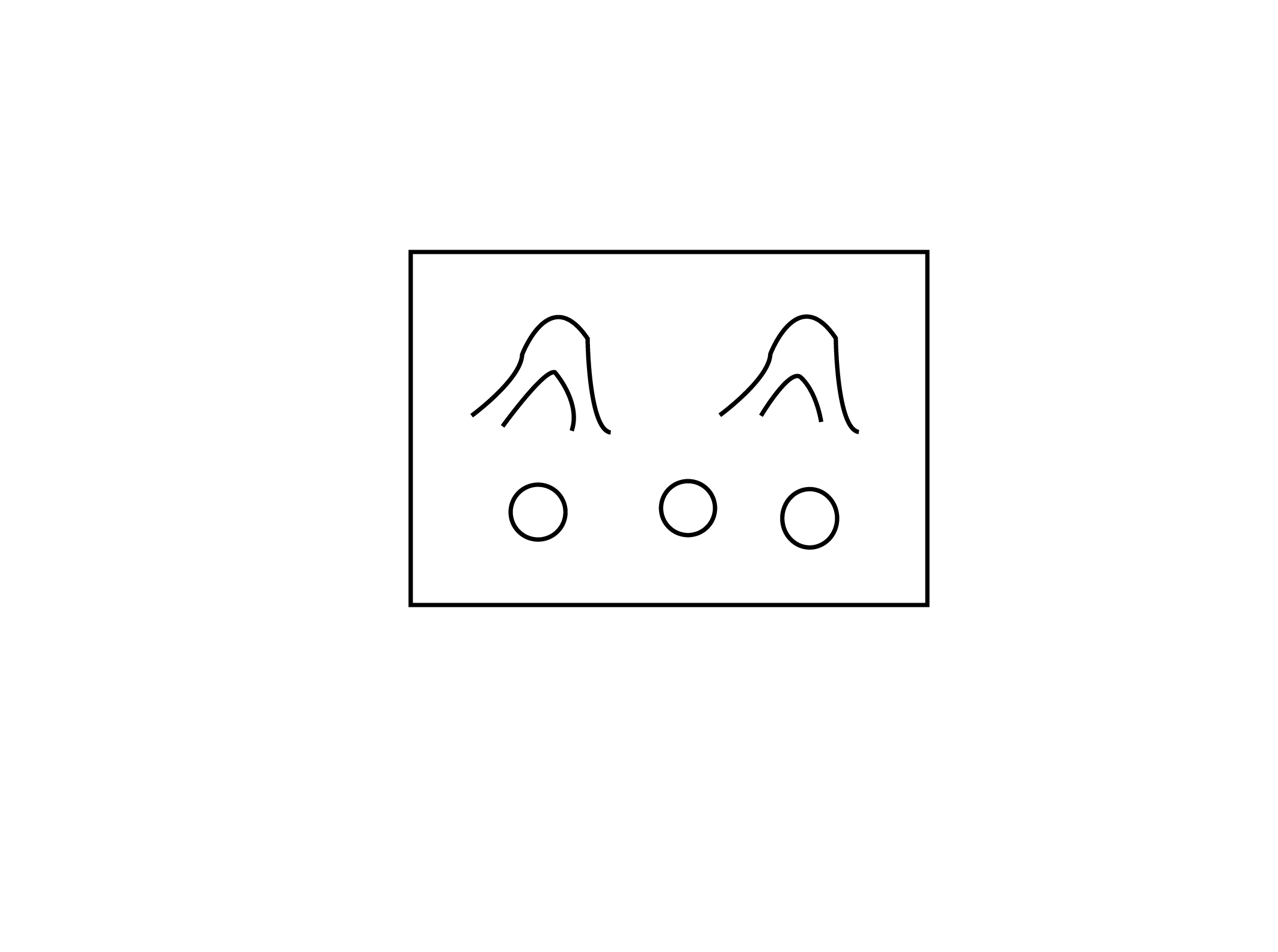}
  \caption{A surface $F \in \cc_{2,4}$ with $\p_0F$ being the rectangular boundary component.}
  \label{fig:figone}
\end{figure}

    Given  an annulus  $C$ as in figure \ref{fig:figtwo} below , then $\sj F$ can be thought of as the boundary connected sum of $F$ with $C$ as pictured in the figure.    

\begin{figure}[ht]
  \centering
  \includegraphics[height=8cm]{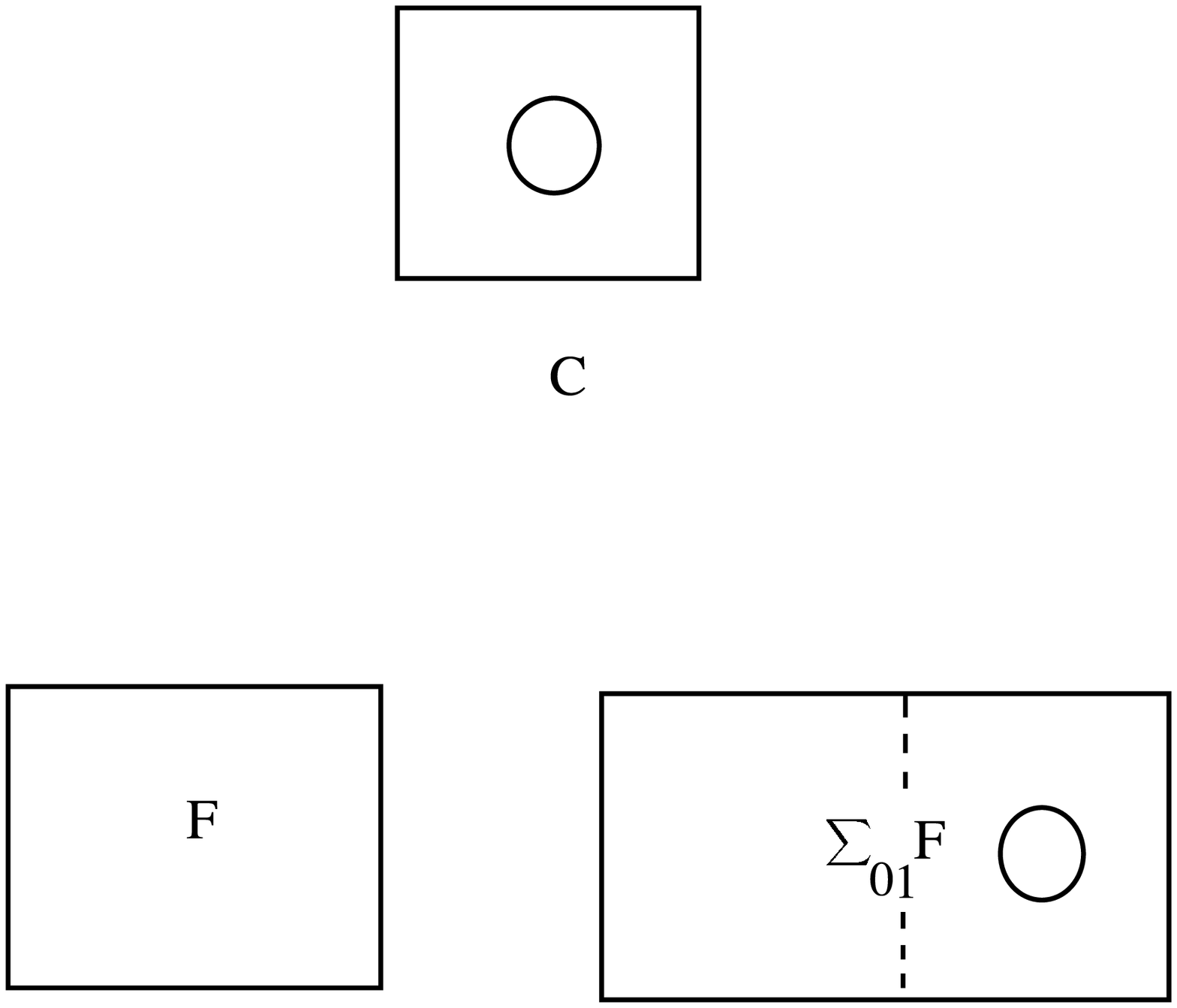}
 \caption{$\sj F$ as a boundary connect sum.  } 
  \label{fig:figtwo}
\end{figure}

The surface $\si F$ can be described similarly, by taking a boundary connect sum with
a surface $D$ of genus one, and one rectangular boundary component, $\p_0D$.   For $F \in \cc_{g,n}$,
the surface $\sk F \in \cc_{g,n-1}$ is obtained by ``filling in the last hole", i.e attaching a disk $D^2$ along
$\p_{n-1}F$ using the parameterization.

   We observe  that the operations $\si$, $\sj$, and $\sk$  are functorial,  since any diffeomorphism $\psi : F_1 \to  F_2$ preserving the boundary parameterizations, induces  diffeomorphisms (preserving boundary parameterizations)  
   
\begin{align}\label{smorph}
  & \si \psi : (\si F_1, \, \phi^{\si F_1}) \to (\si F_2, \, \phi^{\si F_2}),  \\  
   &\sj \psi : (\sj F_1, \phi^{\sj F_1})  \to (\sj F_2, \phi^{\sj F_2}),  \quad \text{and}\notag \\
   &\sk \psi : (\sk F_1, \phi^{\sk F_1})  \to (\sk F_2, \phi^{\sk F_2})
  \end{align}
   defined to be equal to   $e$ on $F_1$,  and the identity on the glued surfaces $C$,  $D$, and $D^2$ respectively.

   We now observe that there are natural embeddings,  yielding morphisms in the category $\cc$, which (by abuse of notation) we also call
   $$
   \si : F \hk \si F   \quad \sj : F \hk \sj F \quad F \hk \sk F
   $$
   These embeddings are  essentially the inclusions 
   of $F$ into the glued surface $\si F$ ,  $ \sj F$, or $\sk F$, together, in the case of $\si$ and $\sj$, with a path from the basepoint in $F$ to the basepoint of $\Sigma F$ which goes below the hole (in the case of $\sj$) or handle (in the case of $\si$).      (See figure  \ref{fig:figthree} below.)

\begin{figure}[ht]
  \centering
  \includegraphics[height=6cm]{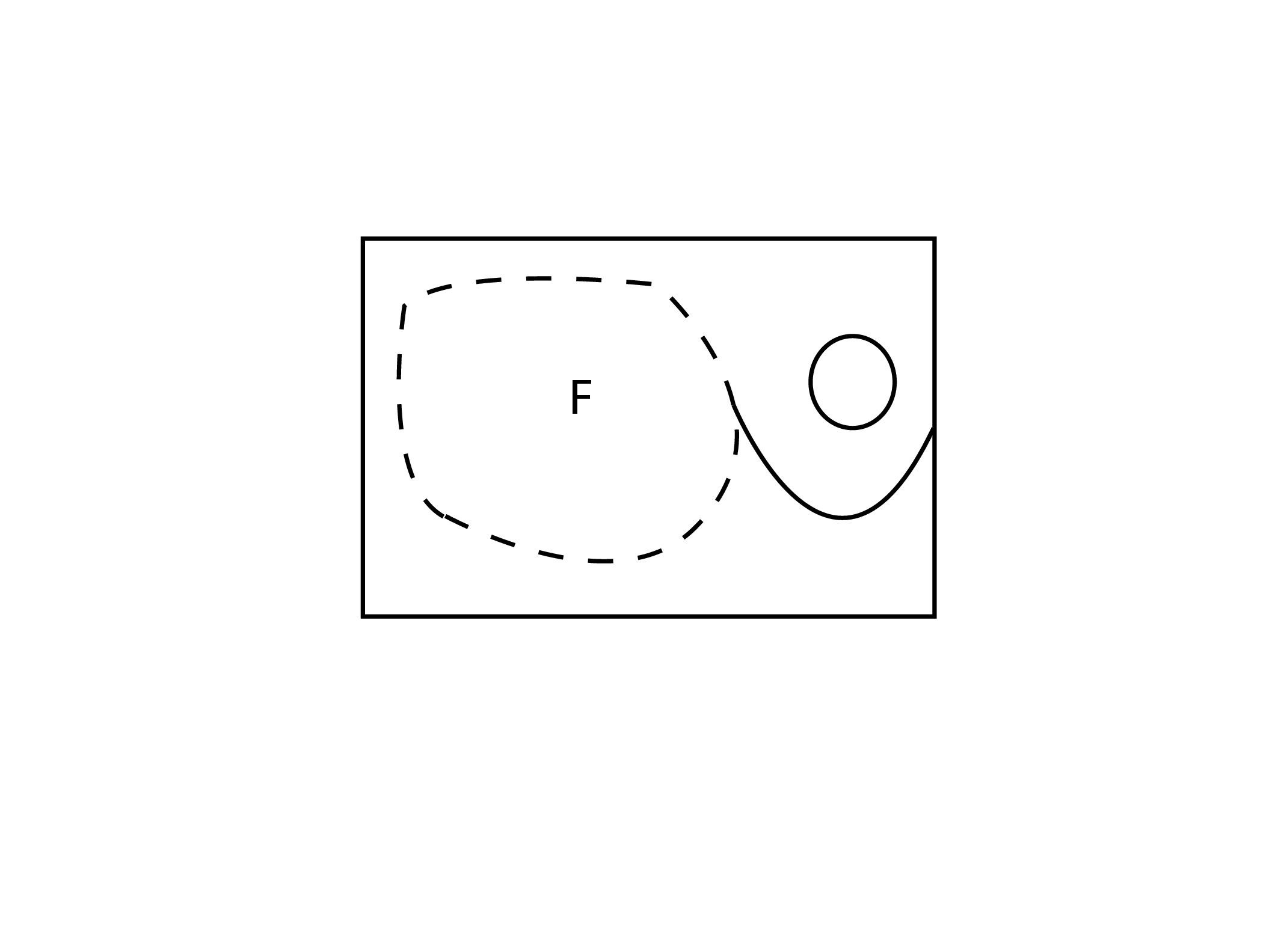}
 \caption{ The morphism $\sj : F \to \sj F$  } 
  \label{fig:figthree}
\end{figure}

       The following is immediate.
   
   \begin{lemma}\label{unique}
   Let $(F, \phi)$ be an object in $\cc_{g,n}$.    Then there are noncanonical  isomorphisms,     
   $$
   F \xr{\cong}\sj \circ \sk F  \quad \text{and} \quad F \xr{\cong} \sk \circ \sj F.
   $$
More generally compositions of the operations $\si$, $\sj$, and $\sk$ give  
isomorphism classes of objects
  $(\Sigma_{i,j} F, \, \Sigma_{i,j}\phi)$ in $\cc_{g+i, j+i}$  for $i \geq 0$ and $j \geq {-n }$.
\end{lemma}
   
   \med
We can actually think of the suspension operations $\si$ and $\sj$ as functors on the entire positive boundary part of the surface category,
   $$
   \si,  \, \sj : \cc_+ \to \cc_+
   $$
   where $\cc_+$ is the full subcategory of $\cc$ generated by surfaces that have nonempty boundaries.  To see how these suspension functors act on morphisms, we use 
   Ivanov's graphical description (see \cite{ivanov}, section 2.5).  Say $e : F_1 \hk F_2$ is a morphism as represented in the left hand of the figure below.  Then $\sj(e)$ is the morphism represented in the right of the figure \ref{fig:figfour}  below.   $\si (e)$ is defined similarly.
   See \cite{ivanov}, section 2.5 for details of this construction.

\begin{figure}[ht]
  \centering
  \includegraphics[height=4cm]{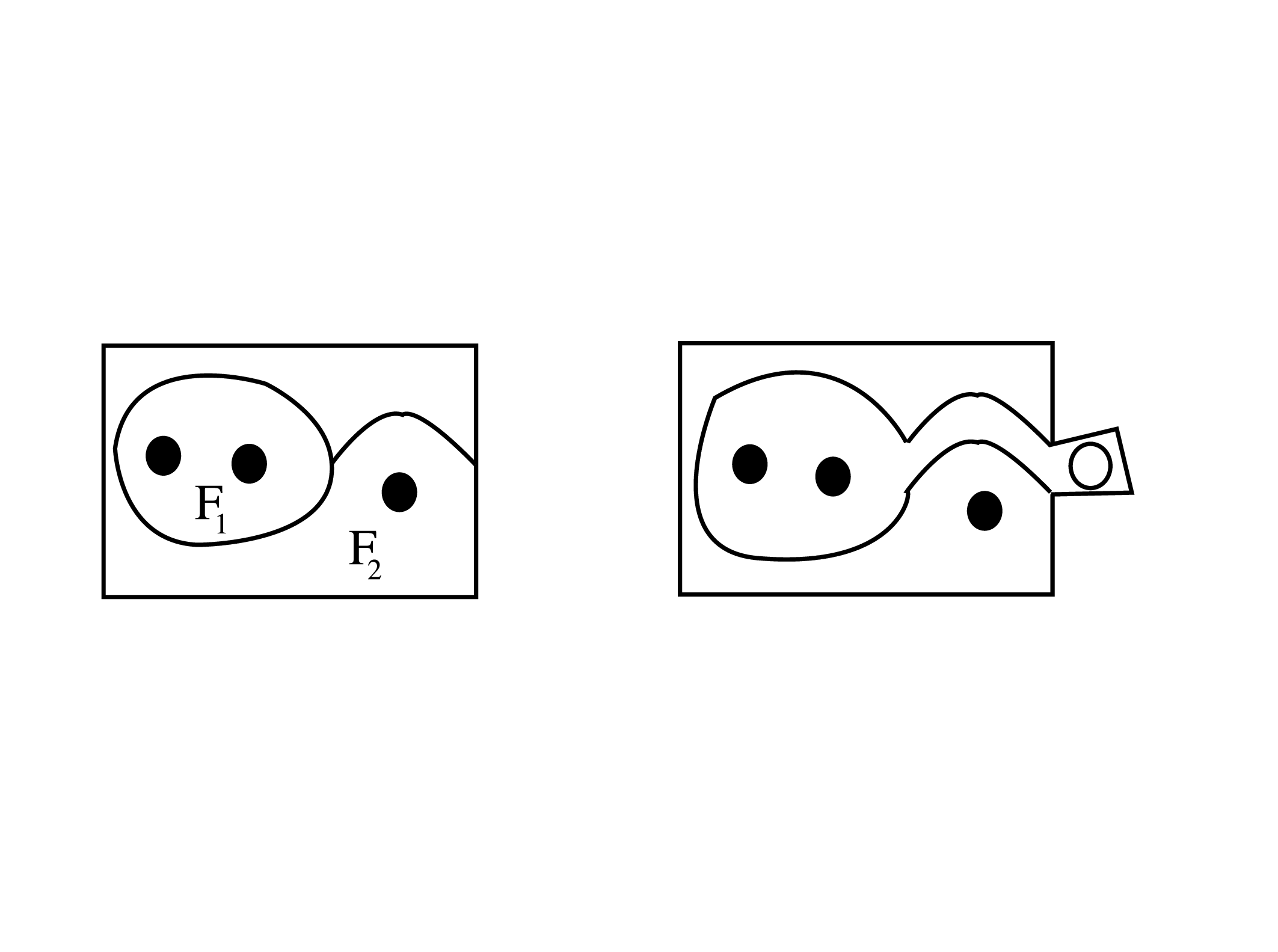}
 \caption{  A morphism $e : F_1 \hk F_2$, and its suspension $\sj (e) : \sj F_1 \hk \sj F_2$.  } 
  \label{fig:figfour}
\end{figure}

   \med
   We now define the notion of a \sl coefficient system. \rm
      
   \begin{definition}\label{coeffsys}
   Let $Ab$ be the category of finitely generated abelian groups and homomorphisms between them.  A    coefficient system $V$  \rm is a covariant functor $V : \cc_+  \to   Ab.$   An extended coefficient system is one that extends to a functor from the entire surface category, $V : \cc \to Ab.$ 
    \end{definition}
   
   
  For an object $(F, \phi)$ of $\cc_{g,n}$ , let $\G (F) = \pi_0(Diff(F, \p F))$ be the mapping class group
  of isotopy classes of orientation preserving diffeomorphisms that fix the boundary pointwise. Notice that in a coefficient system, the $R$-module $ V(F) $ is a module
  over the group ring $\bz [\G
  (F)]$.  (We call it a $\G (F)$-module for short.)
  
  Notice also that given a coefficient system  $V$, one gets new coefficient
  systems $\si V  $ and   $\sj V$,   by defining
  $$
  \Sigma_{i,j} V (F) =   V(\Sigma_{i,j}(F)).
  $$
  On a morphism $e : F_1 \to F_2$,  $\Sigma_{i,j}V(e) : \Sigma_{i,j}F_1 \to \Sigma_{i,j}F_2$ is defined to be $V(\Sigma_{i,j}(e))$, where $\Sigma_{i,j}(e)$ is defined as in figure \ref{fig:figfour}.   Notice that 
  part of the data of a coefficient system consist of natural transformations,
  $$
  \Sigma_{i,j} : V \to \Sigma_{i,j}V.
  $$

  \med
  The following notion of \sl degree, \rm  following Ivanov, who in turn was inspired by the work of Van der Kallen \cite{vanderkallen} and Dwyer \cite{dwyer}, will be very important in our proof of theorem 4.
  
  \med
  \begin{definition}\label{degree}  Let $V$ be a coefficient system. We say that $V$ has degree zero with respect to $\Sigma_{i,j}$ ($(i,j)= (1,0)$ or $(0,1)$)  if 
  the natural transformation, $\Sigma_{i,j} : V(F) \to \Sigma_{i,j}V(F)$ are isomorphisms for all $F$.
  That is,  the  coefficient system is constant.  
  
  Recursively, we say that $V$ has degree $\leq d$ with respect to $\sij$,  if the following two conditions hold:
  \begin{enumerate}
  \item The operation  $\Sigma_{i,j} : V(F) \to \Sigma_{i,j}V(F)$ is a split injection of $\G (F)$-modules, with cokernel $\dij V(F)$.  
  \item The coefficient system $\dij V$ has degree $\leq d-1$ with respect to $\sij$.  
  \end{enumerate}
  
  We say that $V$ has overall degree $\leq d$ if it has degree $\leq d$ with respect to both the functors, $\si$    and $\sj$.
   \end{definition}

   \med
   \bfl 
   \bf Examples.    \rm
   
   \efl
   \begin{enumerate}
   \item  Let $V$ be the coefficient system given as follows.   Let $F \in  \cc_{g,n}$.  Define   $$V(F) = H_1(F) \cong \begin{cases}
  \bz^{2g}  \quad \text{if} \quad n = 0 \\
   \bz^{2g+n-1} \quad \text{if} \quad  n\geq 1. 
   \end{cases}$$  The action of $\G(F)$on $V(F)$ is via the induced map of a diffeomorphism on homology.  
   In this case $\di V(F) = \bz^{2}$ and $\djj V(F) \cong \bz$.  So  both coefficient systems $\di V$and $\djj V$ are constant, and therefore have degree zero.  Therefore  $V$ has degree one  with respect to both $\si$ and $\sj$.  
   
   \item Let $X$ be an Eilenberg-MacLane space, $X = K(A,m)$, with $A$ an abelian  group, and $m \geq 2$.  Defined the coefficient system $V^k$ by
   $$
V^k(F) = H_k (Map_*(F/\p F, X)),
$$

where $Map_*$ denotes the space of based maps.  Here we use the usual conventions
that the basepoint of $F/\p F$ is $\p F$, and that $F/\p F = F_+$, that is, $F$ with a disjoint basepoint, if $\p F = \emptyset$.  We claim that this coefficient system has overall degree $k$. 

To prove this we use induction on $k$.  For  $k = 0$, $$V^0(F) =  H_0(Map_*(F/\p F, K(A,m))) = \begin{cases} \bz[A], \quad \text{if} \quad m=2 \\
\bz \quad \text{if} \quad m>2. \end{cases} $$
In either case  this is a constant coefficient system and so has degree zero.  Now assume inductively the result is true for $V^q$ for $q<k$.  Consider the based homotopy cofibration sequences,
\begin{align}
S^1 \vee S^1 &\to (\si F)/\p (\si F)   \to F/ \p F \notag \\
S^1 &\to (\sj F)/\p (\sj F) \to F/\p F \notag \\
S^0 &\to (\sk F)/\p (\sk F) \to F/\p F \notag 
\end{align}
Here we are using the facts that  $\si F$ is the boundary connect sum of $F$ with the surface $D$ which has the homotopy type of $ S^1\vee S^1$,  and $\sj F$ is the boundary connect sum of $F$ with the surface $C$ which has the homotopy type of $S^1$.   These cofibration sequences induce split fibration sequences,
 \begin{align}
  Map_*(F/\p F, \, X) &\to Map_*((\si F/\p (\si F)), \,  X) \to \Omega X \times \Omega X  \notag \\
  Map_*(F/\p F, \, X) &\to Map_*((\sj F/\p (\sj F)), \,  X) \to \Omega X  \notag \\
   Map_*(F/\p F, \, X) &\to Map_*((\sk F/\p (\sk F)), \,  X) \to X  \notag
 \end{align}
Since $X = K(A,m)$, it has a multiplication, and so the total spaces of these fibrations split  up to homotopy as the products of the fiber and base.  Since $\Omega X$ is  connected, the Kunneth formula gives that $\dij V(F)$ is expressible in terms involving only $H_i(Map_*(F/\p F, \, X)$ for $i \leq k-1, $  as well as  $H_*(X)$ and $H_*(\Omega X)$.  This proves that $V^k$ has degree $\leq k$ with respect to all three functors $\sij$.

We remark that if $X$ is not simply connected, then $H_k(Map_*(F/\p F, \, X) \otimes \tilde H_0(\Omega X \times \Omega X)$ is a direct summand of $\di V^k (F)$, so $V^k$ has  infinite degree with respect to $\si$.  It also  has infinite degree with respect to $\sj$.  Indeed this is the primary reason that we need $X$ to be simply connected in the statement of Theorem \ref{one}.
   
   \end{enumerate}

   Given a coefficient system,  Theorem  \ref{mpg}  is about the homology groups $H_q(\G (F); V(F))$. 
  Following the notation of Ivanov \cite{ivanov},  we define the relative homology group,
   $$
Rel^V_q (\sij F, F) = H_q((\G(\sij F), V(\sij F), (\G(F), V(F)).
$$

\med
\bfl
\bf Remark.  \rm Even though the functor $\sk : \cc_{g, n} \to \cc_{g, n-1}$ does not extend to a functor on all of $\cc$,   or even on $\cc^+$, we can still define these relative homology groups with respect to $\sk$, for the following reason.  Given a coefficient system,
$V : \cc \to Ab$, one still has a natural transformation  $\sk$ between the restriction of $V$ to $\cc_{g,n}$ and the composition $\cc_{g,n} \xr{\sk} \cc_{g,n-1} \xr{V} Ab$ for each $(g,n)$ with $n>0$.  This allows us to define the relative homology groups,
$$Rel_q^V(\sk F, F) = H_q(\Gamma (\sk F), V(\sk F);  \, \Gamma (F),  V(F))  $$ for any surface $F \in \cc_{g,n}, \, n>0$. 

More generally, for any surface $F \in \cc_{g,n}, \, n>0,$ we can define relative homology groups in the following way.  Consider sequences of pairs, $I = ((i_1,j_1), \cdots , (i_k,j_k))$, with each $(i_r, j_r)$ of the form $(1,0), \, (0,1),$ or $(0, -1)$.  We call such a sequence \sl admissible \rm if $j_1 \leq j_2 \leq \cdots \leq j_k$, and if $ j = \sum_{r=1}^k j_r$, then $n+j > 0$.  (By convention, the empty sequence is also called admissible.) We let $\Sigma^I F = \Sigma_{(i_1, j_1)} \circ \cdots \circ \Sigma_{(i_k, j_k)}(F)$.  Using the above mentioned  natural transformations,  we may define the relative homology groups,
$$
Rel^{\Sigma^IV}_q (\sij F, F) = H_q((\G(\sij \Sigma^I F), V(\sij \Sigma^I F); \, (\G(\Sigma^I F), V(\Sigma^I F)).
$$

\efl
Consider the  following long exact sequence. 

$$
\cdots \to H_q(\G (F); V(F)) \to H_q(\G (\sij F); V(\sij (F)) \to \rv_q(\sij F, F) \to H_{q-1}(\G(F); V(F)) \to \cdots
$$
 From this sequence Theorem \ref{mpg}    is immediately seen to be a consequence of the following stability theorem.

\med
\begin{theorem}\label{mpg2}
Let $V$ be a coefficient system of overall degree $\leq d$.  Let   $F \in \cc_{g,n}, \, n>0$.  Then for any admissible sequence $I$, the relative groups,
$$Rel^{\Sigma^IV}_q (\sij F, F)= 0$$ for  
for $2q+d+2 \leq g$, when  $(i,j) = (1,0), (0,1),$  and for $2q+d+2 \leq g-1$ when $(i,j)=(0,-1)$.  In this last case ($(i,j)=(0,-1)$) $V$ is assumed to be an extended coefficient system, so that it is defined for closed surfaces, as well as surfaces with boundary. 
\end{theorem}

\med
We remark  that for $(i,j) = (1,0)$, this theorem was proved by Ivanov \cite{ivanov}.\footnote{ Ivanov assumed that the surface $F$ had only one boundary circle.  However as he pointed out at the end of his paper, his argument easily extends to the general case.}
Our proof of this more general theorem follows the ideas of Ivanov, but in the cases of
$(i,j) = (0,1)$ and $(0,-1)$ it will require  further argument.

Our proof of this theorem goes by induction on the degree $d$ of the coefficient system $V$.
Degree zero coefficient systems are constant, with trivial mapping class group action.  In this case
 Theorem  \ref{mpg2}, and in particular,  Theorem  \ref{mpg} is the stability theorem of Harer \cite{harer}
as improved by Ivanov \cite{ivanov}.  

\med

In what follows, we write $Rel^V_q (\sij F, F)$ to denote any of the relative groups, $Rel^{\Sigma^IV}_q (\sij F, F) $ for any admissible sequence $I$ (including the empty sequence).

\med
\begin{induct}\label{assume1}
We inductively assume  Theorem  \ref{mpg2} to be true for coefficient systems of overall degree $< d$.
\end{induct}

\med
Our strategy for the completion of this inductive step, and thereby the completion of the proof of  Theorem  \ref{mpg2} is the following.  As mentioned above, we already know this theorem to be true for $\Sigma_{1,0}$ by the work of Ivanov \cite{ivanov}.  We will then complete the inductive step for the natural transformations $\sj$ and $\sk$ separately.  

In both cases, our arguments will rely on the action of the relevant mapping class groups
on certain simplicial complexes (the ``curve complexes" of Harer \cite{harer}) that are highly
connected.  We will then analyze the corresponding spectral sequence.   However this analysis
(and indeed the choice of curve complexes) is a bit different  in the cases of the two transformations
$\sj$ and $\sk$, which is why we deal with them separately.

In sections 1.2 and 1.3 we complete the inductive step for the operation $\sj$.  We begin section 1.2 by recalling the spectral sequence of a group action on a simplicial complex, and then describe
the curve complex we will study.  We analyze the spectral sequence and the upshot
is a result giving   stability of the relative homology groups $\rv_q(\sj F, F)$.  Of course we want to prove that these groups are zero, and this requires further argument, which is done in section 1.3.
In section 1.4 we complete the inductive step for the operation $\sk$, using an action on a slightly different curve complex. 

Throughout the rest of chapter 1 we will be operating under Inductive Assumption \ref{assume1}.
 
\subsection{The curve complex and a relative spectral sequence}

In this section and section 1.3 our goal is to complete the inductive step in Assumption \ref{assume1}, for the operation $\sj$.  That
 is, we want  to prove that for a coefficient system $V$ of degree $ \leq d$ then the relative groups $\rv_q(\sj F,F)  = 0$ for $2q \leq g-d$.   To do this, we use another argument, this time
 induction on $q$.  We assume the following.
 
 \begin{induct}\label{assume2}
 Let $V$ be a coefficient system of degree $\leq d$. Then for $q < m$,   $\rv_q(\sj F,F)  = 0$ for any surface $F$ of genus $g\geq 2q+d+2$. 
 \end{induct}
 
 Clearly by completing this inductive step, we will complete the inductive step for assumption \ref{assume1}, and thereby complete proof of  Theorem  \ref{mpg2} for the operation $\sj$.
 
So in the following two sections we will operate under Inductive Assumption \ref{assume2} within 
assumption \ref{assume1}.  Our goal is to prove that $\rv_m(\sj F, F) = 0$.  We do this in two steps. The first step, which is the object of this section, is to prove the following.

\begin{proposition}\label{isom} Let $F$ be any surface with boundary of genus $g \geq 2m+d+2$.    Then  there is an isomorphism 
$$
\rv_m( (\sj F , F) \xr{\cong} \rv_m(\si^2F, \si F). 
$$
    
\end{proposition}

The second step, which we complete in section 1.3, is to prove that under these hypotheses,  $\rv_m(\sj F, F) = 0$.

\subsubsection{The spectral sequence of a group action}
We begin our proof of Proposition \ref{isom} by recalling the spectral sequence for an action of a discrete group $G$ on a simplicial complex.  See \cite{brown} or  section 1.4 of  \cite{ivanov} for a more complete description.

Let $X$ be a simplicial complex with an action of $G$.  In particular $G$ acts on the set of $p$-simplices, for each $p$.  Given a simplex $\sigma$, let $G_\sigma$ denote the stabilizer subgroup of this action.  We let $\bz_\sigma$ be the ``orientation $G_\sigma$-module", defined to be $\bz$ additively, with the action of $g \in \G_\sigma$ to be multiplication by  $\pm 1$ depending on whether $g$ preserves or reverses orientation.

For $M$ a $G$-module, let $M_\sigma = M \otimes \bz_\sigma$.   For each $p$, let $Simp_p$ be a set of representatives of the orbits of the action of $G$ on the $p$-simplices of $X$.  Then applying equivariant homology to the skeletal filtration of $X$ defines a spectral sequence whose $E^1$ term is given by

\begin{equation}\label{spectral}
E^1_{p,q}  = \bigoplus_{\sigma \in Simp_{p}} H_q(G_\sigma, M_\sigma)
\end{equation}
which converges to the equivariant homology, $H_{p+q}^G(X; M).$    

\med
There is a relative version of this spectral sequence as well.   This applies when considering two spaces:  a $G$-simplicial complex $X$, and a $G'$-simplicial complex $X'$.  Suppose $\phi : G \to G'$ is a homomorphism, and $f : X \to X'$ is a simplicial, equivariant map with respect to the homomorphism $\phi$.  Suppose also
that $\psi : M \to M'$ is an $\phi$-equivariant homomorphism between the $G$-module $M$ and the $G'$-module $M'$.  Then using
the mapping cylinder, one can define the relative groups, $H_*^{G', G}(X', X; M', M)$.

Now assume that the map $ f : X \to X'$ has the special property that it induces a bijection on the $G$-orbits of the $k$-simplices of $X$ to the $G'$-orbits of the $k$-simplices of $X'$, for each $k$.   If $X$ and $X'$ are $d$-connected, then   there is a relative spectral sequence whose $E^1$-term is given by

\begin{equation}\label{relspec}
E^1_{p,q} = \bigoplus_{\sigma \in Simp_{p-1}} H_q(G'_\sigma, G_\sigma; \, M', M)  
\end{equation}
which converges to  zero in the range $p+q \leq d+1$.  Here we formally  let $Simp_{-1}$ consist of one element, $\sigma_{-1}$, and let $G_{\sigma_{-1}} = G$, $G'_{\sigma_{-1}} = G'$.  See \cite{ivanov},  \cite{vanderkallen} for details of this spectral sequence.  

We will use this spectral sequence to prove Proposition $\ref{isom}$.

\subsubsection{The curve complex}
Let $(F, \phi)$ be a fixed  object in $\cc_{g, n}$, with $n \geq 2.$  We recall the definition of the the curve complex  $H(F)$ studied by Harer \cite{harer} and Ivanov \cite{ivanov}.  This will be a $\G_{g,n}$-equivariant simplicial complex, and we shall apply the relative spectral sequence (\ref{relspec}).   For ease
of notation let $C_0 =\p_0 F$, and $C_1 = \p_1 F$.   Choose fixed points $b_0 \in C_0$, and $b_1 \in C_1$.  

The simplicial complex $H(F) = H(F; b_0, b_1)$ is the complex whose vertices are represented by isotopy classes of embedded  arcs in F    from $b_0$ to $b_1$.  The interior of these arcs must lie in the interior of $F$.  A set $\{ \alpha_0, \cdots , \alpha_p\}$ of such vertices spans a $p$-simplex in $H(F)$ if there are disjoint representing arcs,
$A_0, \cdots , A_p$ such that $F- \cup_i A_i$ is connected.   The dimension of $H(F)$ is $2g$, and  it was shown in \cite{harer}, \cite{wahl} that this space is  $(2g-1)$-connected. Therefore it is   homotopy equivalent to a wedge of $2g$-dimensional spheres.  The mapping class group $\G (F)$ acts on $H(F)$ since it acts on the set of vertices, and maps simplices to simplices.   However the action is not transitive on the set of $p$-simplices, if $p \geq 1$.  Indeed the orbit of a given $p$-simplex $A = (A_0, \cdots, A_{p})$ is determined by a certain permutation $\sigma (A) \in \Sigma_{p+1}$.  Moreover if $g$ is sufficiently large 
 all such permutations can appear in this way. 

The permutation $\sigma (A)$ is defined as follows.   We orient the two boundary circles $C_0$ and $C_1$ so that $C_0$ is incoming and $C_1$ is outgoing.  Then the trivialization of the normal bundle of an arc from $b_0$ to $b_1$ determined by $C_0$ at $b_0$ agrees with the trivialization  determined by $C_1$ at $b_1$.    Given a $ p$-simplex $A=(A_0, \cdots, A_{p})$ of arcs  that start  at $b_0$ in the order given by the orientation of $C_0$,  they arrive at $b_1$ in an order which is the permutation $\sigma (A)$ of the order at $b_1$ dictated by the orientation of $C_1$.  See figure \ref{fig:figfive} below.

\begin{figure}[ht]
  \centering
  \includegraphics[height=4cm]{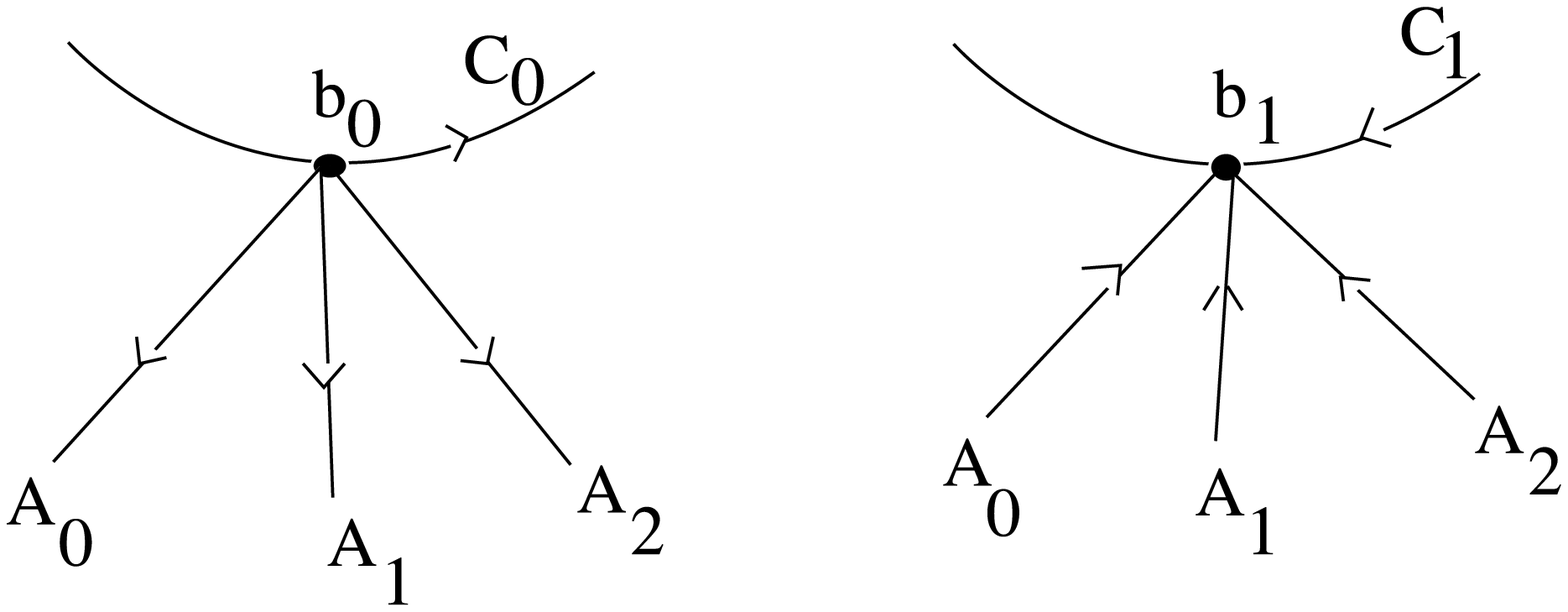}
 \caption{$\sigma (A) = (0,2) \in \Sigma_3$ } 
  \label{fig:figfive}
\end{figure}
 To see that $\sigma (A)$ determines the orbit type, notice that it determines the boundary structure of the cut surface, $F-A$.   The boundary of $F-A$ consists of $\p F - (C_0 \cup C_1)$ together with the new boundaries created by the cuts.  (See the example below.)  If $\sigma (A) = \sigma (A')$, then  $F-A$ and $F-A'$ have the same new boundary structure, and there is a preferred diffeomorphism between their boundaries.  This extends to a diffeomorphism between $F-A$ and $F-A'$ since they have the same Euler characteristic ($\chi (F-A) = \chi (F-A') = \chi (F) +p+1$).  Sewing the cuts back together   defines an element $\gamma \in \G (F)$ with $\gamma (A) = A'$.  
 
 \med
 \bf Examples.  \rm   Suppose $F$ has genus $g$ and $r+2$ boundary components, and let $A = (A_0, \cdots , A_p)$ be a $p$-simplex of $H(F)$.   
 \begin{enumerate}\label{ident}
 \item Suppose $\sigma (A) = id.$  The new boundaries are
 $(A_0, \bar C_1, \bar A_p, C_0), (\bar A_0,  A_1), \cdots, (\bar A_{p-1},  A_p)$, and $F-A$ has genus $g-p $ and $r+p-1$ boundary components.   In this notation, and in what follows, given an oriented edge $E$, $\bar E$ denotes the edge with the opposite orientation.
 \item  Suppose $\sigma (A) =\prod_{i=0}^{[p/2]} (i, p-i) \in \Sigma_{p+1}$.  The new boundaries are: $$\begin{cases} (A_0, \bar A_1, \cdots , A_p, C_0) \, \text{and} \, (A_0, C_1, \bar A_p, \cdots,  A_1) \quad \text{if $p$ is odd,} \\
 (A_0, \bar A_1, \cdots. A_p, \bar C_1, \bar A_0, \cdots, A_p, C_0) \quad \text{if $p$ is even .}  \end{cases}
 $$ In this case $F-A$ has genus  $g - \frac{p}{2}$  with $r+1$ boundary components when $p$ is even, and genus $g-\frac{p-1}{2}$ with $r+2$ boundary components when $p$ is odd.    
 \item In general the number of boundary components in $F-A$ is $r+t +1$, with $0\leq t \leq p$.  Since the Euler characteristic $\chi (F-A)=\chi (F) +p+1$, one gets that
 $2g(F-A)=2g-p-t$, and hence $g(F-A)\geq g-p$.  
  \end{enumerate}

 \med
 Let $F$ have genus $g$ and $r$ boundary components, with $r \geq 1$.    Let $R = \sj F$.   Let $H(\sj F)$ be the arc complex with respect to $b_0, b_1 \in \p (\sj F)$. (As above, $b_0 \in \p_0 (\sj F)$, and $b_1 \in \p_1(\sj F)$.)   Choose $b_1' \in \p_1(\sj R)$, and choose an embedded path from $b_1$ to $b_1'$ in $ \sj R - interior (R)$.  See figure \ref{fig:figsix} below.
 
\begin{figure}[ht]
  \centering
  \includegraphics[height=6cm]{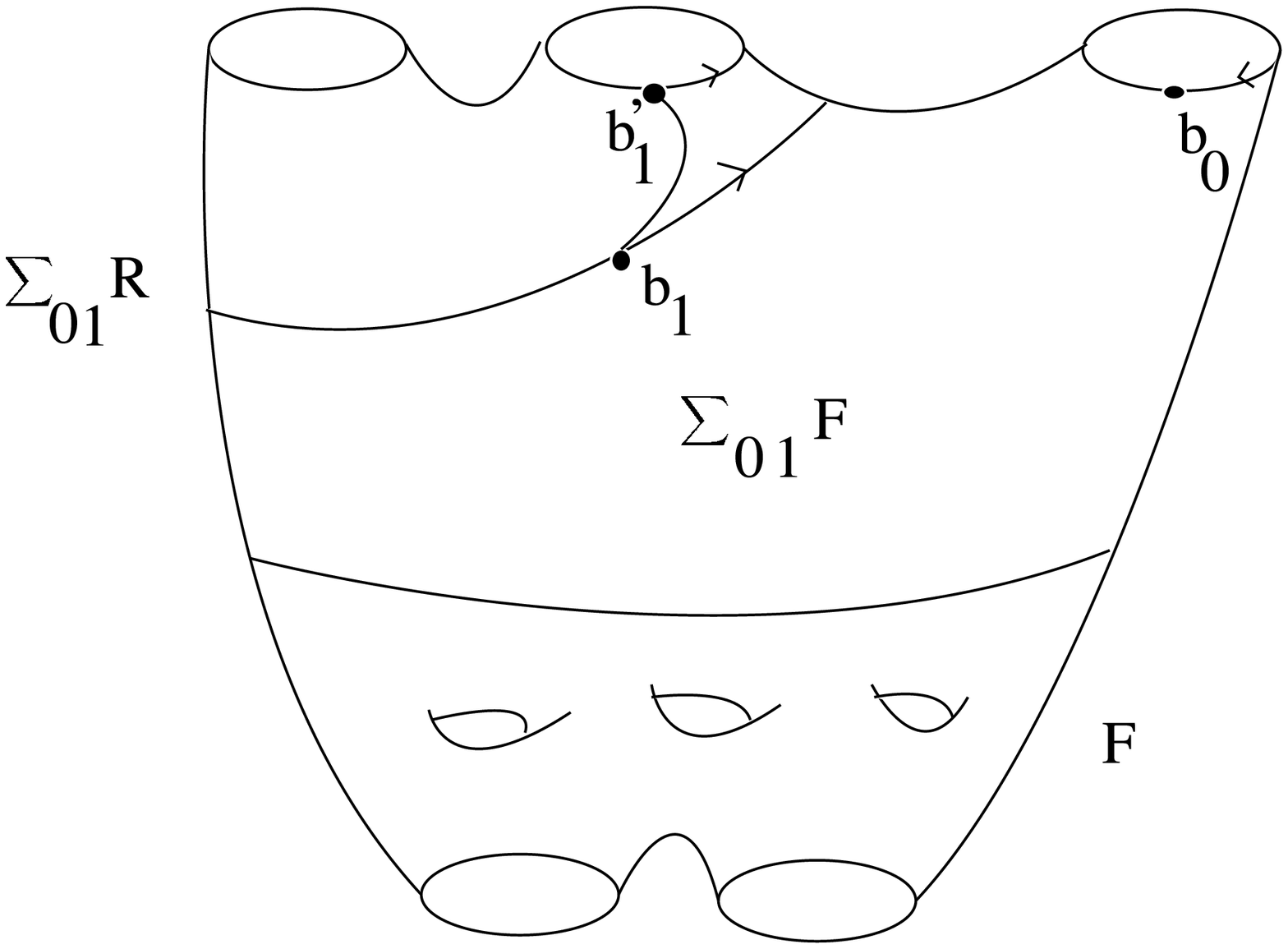}
 \caption{$\sj R$ } 
  \label{fig:figsix}
\end{figure}

The path (or rather a closed tubular neighborhood around the path) induces a simplicial map
$$
\iota : H(\sj F) \la H(\sj R)
$$
which is equivariant with respect to 
   $$
   \iota : \G(\sj F) \la \G(\sj R),
   $$ and such that the permutation $\sigma (A)
   $ associated with a $p$-simplex of $H(\sj F)$ is equal to the permutation associated with $\iota (A)$ in $H(\sj R)$.  In particular $\sj R - \iota (A)$ has precisely one more boundary component that $\sj F - A$, and the two cut surfaces have the same genus.  The examples above show  that  the genus of $\sj F - A$ is $\geq g-p+1$, and the number of boundary components is $r +t$, where $0\leq t \leq p-1$.  In particular, $\iota$ preserves the number of orbits of $p$-simplices.  
   
   \med
   We now consider the 
    relative spectral sequence (\ref{relspec}) as it pertains to this situation.   This spectral sequence  has $E^1$-term equal to 
   \begin{equation}\label{e1}
   E^1_{p,q} = \bigoplus H_q(\G (\sj R)_{\iota (A)}, \, \G (\sj F)_A; \, V(\sj R), \, V(\sj F)),
   \end{equation}
   where  the direct sum varies over the $\G (\sj F)$-orbits of $(p-1)$-simplices $A \in H(\sj F)$.   Note that
   
   \begin{align}
   E^1_{0,q} &= H_q(\G(\sj R), \, \G(\sj F); \, V(\sj R), \, V(\sj F)) = \rv_q(\sj R, \sj F) \notag \\
   E^1_{1,q} &= H_q(\G (R), \, \G(F); \,  \, V(\sj R), \, V(\sj F)) = Rel^{\Sigma_{0,1}V}_q(R, F). \notag
   \end{align}
   
   We are interested in the spectral sequence in total degrees $\leq 2g$.  In this range $E^\infty_{p,q} = 0$
   since $H(\sj F)$ and $H(\sj R)$ are both $(2g-1)$-connected.  
   
   The  isotropy groups that appear in (\ref{e1}) are isomorphic to the mapping class groups of the cut surfaces, $\sj F-A$, and 
   $\sj R - \iota(A)$.    Our goal in proving Proposition \ref{isom}  is to construct a type of suspension map 
   $$
   \rv_m(R,F) \to \rv_m(\sj R, \sj F)
   $$
 which  is an isomorphism.  Here $m$ is fixed by Inductive Assumption \ref{assume2}.  Recall this in particular implies that $g \geq 2m+d+2$.  
   This map  will be  the composition of two maps, namely
   
   $$
    \rv_m(R,F) = H_m(\G(R), \G (F); \, V(R), V(F)) \to H_m (\G(R), \G(F); \, V(\sj R), V(\sj F))
    $$
    and the differential,
    $$
    d^1 : E^1_{1,m} \to E^1_{0,m}
    $$
    in the spectral sequence.     The first of these two maps is injective, with cokernel $$ H_m(\G(R), \G (F); \, \djj V(R), \djj V(F)).$$   But by Inductive Assumption \ref{assume1}, this group is zero because $\djj V$ is a coefficient system of degree $\leq d-1$.  Thus to prove Proposition \ref{isom} we need to show the differential $d^1: E^1_{1,m} \to E^1_{0,m}$ is an isomorphism.

    \med
    Now the $E^1$-term in the spectral sequence is a sum of the relative homology groups,
    $$
    \bigoplus_{A  \in Simp_{p-1}} H_q(\G (\sj R -  \iota (A)), \, \G (\sj F - A); \, V(\sj R), \, V(\sj F)).  
    $$

Assume $p \geq 2$ and  $p+q=m+1$. For a fixed $(p-1)$-simplex $A$, let $g' = genus (\sj F-A)$.  Notice that  $g \geq g' \geq g-p+1$.  Now    $q =  m+1-p  \leq  \frac{g-d}{2} + 1-p.$   Since $g-p+1 \leq g'$, $\frac{g-d}{2} + 1-p \leq g' - \frac{g+d}{2}.$  So $q \leq g' - \frac{g+d}{2} \leq \frac{g'-d}{2}.$ (This last inequality follows from the fact that $g' \leq g$.)   Thus the Inductive Assumption \ref{assume2} implies that for $p \geq 2$,
$E^1_{p,q} = 0$ for $p+q =m+1$.   Since $E^\infty_{p,q} = 0$ for $p+q = m$, this implies that the differential

   $$
   d^1 : E^1_{1,m} \la E^1_{0,m}
   $$
   must be surjective. If not, the cokernel of $d^1$ would survive to $E^\infty_{0,m}$.     Hence
   $$
   \rv_m(R,F) \la \rv_m(\sj R, \sj F)
   $$
   is surjective.  But it is also injective, since there is a right inverse induced by $\sk$,
   $$
    \rv_m(\sj R, \sj F) \la    \rv_m(R,F).
    $$
    Recalling that $R = \sj F$, we have now proved Proposition  \ref{isom}.

   \subsection{The completion of the inductive step}   
   In this section our goal is to complete the inductive step and thereby complete the proof of  Theorem  \ref{mpg2} for the operation $\sj$.    We continue to operate under Inductive Assumptions \ref{assume1} and \ref{assume2}.   We need to prove the following. 
   
  \begin{proposition}\label{complete} Let $V$ be a coefficient system of degree $\leq d$.  Let $F$ be any surface with boundary of genus $g$ with $g\geq 2m+d+2$.  Then $\rv_m(\sj F, F) = 0$.  
  \end{proposition} 

\med
Before we begin the proof we need some preliminary results.

\med
We begin by defining embeddings,  $d_i : \sj F \hk \sj^k F$, for $i = 0, \cdots, k$. 
To describe these embeddings  it is helpful for graphical reasons  to continue to think of the suspension functor  $\sj$ as described by figure \ref{fig:figthree} in section 1.1 above.   This way of depicting the functor $\sj : \cc_{g,n} \to \cc_{g, n+1}$ allows us to describe
embeddings, $d_i : \sj F \hk \sj^k F$ for $i = 0, \cdots k$.  These embeddings are described by figure \ref{fig:figseven} below:

\begin{figure}[ht]
  \centering
  \includegraphics[height=10cm]{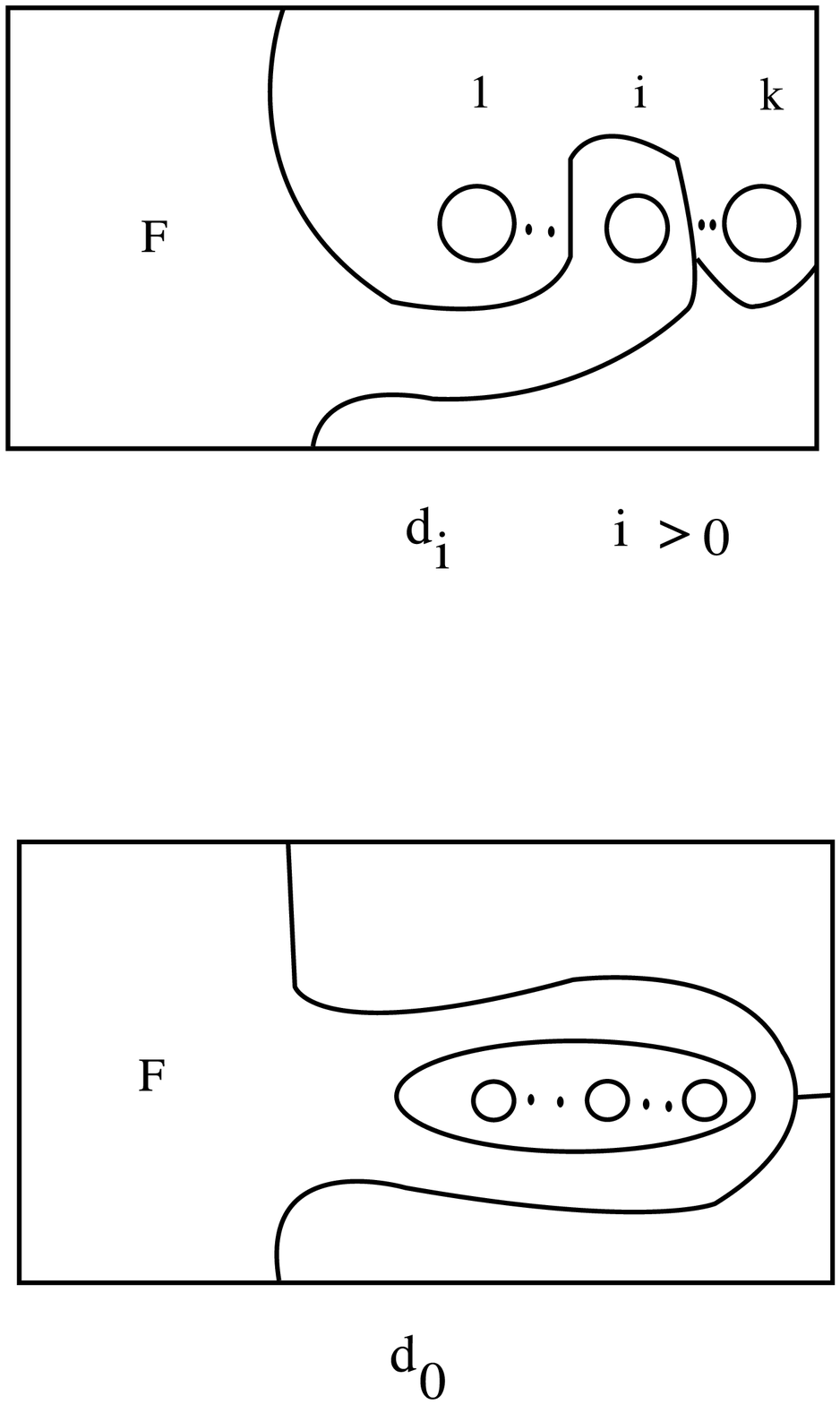}
 \caption{The embeddings $d_i: \sj F \to \sj^k F$, $i = 0, \cdots , k$.  } 
  \label{fig:figseven}
\end{figure}

We view  each map   $d_i:  \sj F \hk  \sj^k F$ as a morphism in the surface category $\cc$.    They therefore define homomorphisms of the relative homology groups,
$$
d_i : \rv_m(\sj F, F) \to \rv_m(\sj^k F, F).
$$

Now for $i = 1, \cdots , k$,  let $e_i : \sj^k F \to \sj F $ be the embedding that attaches a disk to all but the $i^{th}$ hole in in $\sj^k F$.   These maps are also  viewed as morphisms in the surface category $\cc$, and notice that $e_i \circ d_i  : \sj F \to \sj F$ is isotopic to the identity morphism.  We also observe that by figure 7 below, the  following is clear.

\begin{proposition}\label{relations}
   \begin{enumerate}
   \item For $i = 1,\cdots, k$, then   the homomorphism 
 $$
d_i : \rv_m(\sj F, F) \to \rv_m(\sj^k F,   F),
$$
is a  split injection  with left inverse $e_i :  
\rv_m(\sj^k F, \sj F) \to \rv_m(\sj F, F)$. 
\item  For $i = 1, \cdots, k$ and $  j\neq i,$  the composition
$$
e_j \circ d_i : \rv_m(\sj F, F) \to \rv_m(\sj F, F)
$$
is the zero homomorphism.  
\item  For $i = 0$, and $j = 1, \cdots k$,
$$e_j \circ d_0  : \rv_m(\sj F, F) \to \rv_m(\sj F, F)$$ is the identity.
\end{enumerate}
\end{proposition}

\med
\begin{lemma}\label{sumdees}
For each $k$,  the map
$$
\oplus_{i=1}^k d_i : \bigoplus_{i=1}^k \rv_m(\sj F, F) \la  \rv_m(\sj^k F, F)
$$
is an isomorphism, with inverse, $\oplus_{i=1}^k e_i$.  
\end{lemma}

\med
\begin{proof}  We begin by observing that in the case $k = 2$, Proposition \ref{isom}
says that the composition
$$
\bar d_1 : \rv_m(\sj F, F) \xr{d_1} \rv_m(\sj^2 F, F) \xr{\pi_2} \rv_m(\sj^2 F, \sj F)
$$
is an isomorphism, where $\pi_2$ is the projection in the short exact sequence,
$$
0 \la \rv_m(\sj F, F) \xr{d_2} \rv_m(\sj^2 F, F) \xr{\pi_2} \rv (\sj^2 F, \sj F) \la 0.
$$
This implies that $$d_1 \oplus d_2 : \rv_m(\sj F, F) \oplus \rv_m(\sj F, F) \to \rv_m(\sj^2 F, F)$$ is an isomorphism, as stated in the lemma. 

\med
For general $k$, the argument is similar, but somewhat more complicated.   First notice that by
 Proposition \ref{relations} part (1) there is a split short exact sequence,
\begin{equation}\label{deeone}
0 \to \rv_m(\sj F, F) \xr{d_1}  \rv_m(\sj^k F, F) \to \rv_m(\sj^kF, \sj F) \to 0,
\end{equation}
and an identification of $ker (e_1) \cong \rv_m(\sj^kF, \sj F)$.  Since, by part (2) of the proposition, $e_1 \circ d_2 = 0$, we have that 
$$
d_2 : \rv_m(\sj F, F) \to ker (e_1) \cong \rv_m(\sj^kF, \sj F) \hk  \rv_m(\sj^k F, F)
$$
is split injective, split by $e_2 : ker (e_1) \subset  \rv_m(\sj^k F, F) \to  \rv_m(\sj F, F)$.
Using Proposition \ref{isom}, we have an identification of  $\rv_m(\sj F, F)$ with
$ \rv_m(\sj^2 F, \sj F)$, and so we get a split short exact sequence,
\begin{align}\label{deetwo}
0 \to  \rv_m(\sj F, F) \cong  \rv_m(\sj^2 F, \sj F) &\xr{d_2}  \rv_m(\sj^kF, \sj F)\cong ker (e_1) \\
&\to  \rv_m(\sj^kF, \sj^2 F)\to 0 \notag
\end{align}
and an identification of $ker (e_1) \cap ker (e_2) \cong \rv_m(\sj^kF, \sj^2 F)$.  Putting equations (\ref{deeone}) and (\ref{deetwo}) together, we have a split short exact sequence,
 
$$
0 \to \bigoplus_{i=1}^{2} \rv_m(\sj F, F) \xr{d_1 \oplus d_2} \rv_m(\sj^k F, F) \to \rv_m(\sj^k F, \sj^{2} F) \to 0.
$$

Continuing in this way, we have, for each $j$,  a  split short exact sequence

$$
0 \to \bigoplus_{i=1}^j \rv_m(\sj F, F) \xr{\oplus_{i=1}^j d_i} \rv_m(\sj^k F, F) \to \rv_m(\sj^k F, \sj^j F) \to 0
$$
and an identification of $\bigcap_{i=1}^j ker (e_i)   \cong \rv_m(\sj^k F, \sj^j F)$. 

At the final stage we have a split short exact sequence
$$
0 \to \bigoplus_{i=1}^{k-1} \rv_m(\sj F, F) \xr{\oplus_{i=1}^{k-1} d_i} \rv_m(\sj^k F, F) \to \rv_m(\sj^k F, \sj^{k-1} F) \to 0
$$
and a split injective map $d_k : \rv_m(\sj F, F) \to \bigcap_{i=1}^{k-1} ker (e_i)  \cong \rv_m(\sj^k F, \sj^{k-1} F).$   But by Proposition \ref{isom} these two groups are isomorphic.  Thus $d_k : \rv_m(\sj F, F) \to \rv_m(\sj^k F, \sj^{k-1} F)$ is an isomorphism,
and the lemma is proved. 
\end{proof}

\med
This lemma allows us to prove the following, which we will use later in the argument.

\med
\begin{corollary}\label{moredees}
\begin{enumerate}
\item The following diagram commutes:
$$
\begin{CD}
\rv_m(\sj F, F)  @>d_0 >> \rv(\sj^k F, F) \\
@V= VV   @VV \oplus_{i=1}^k e_i V \\
\rv_m(\sj F, F) @>>\Delta >  \bigoplus_{i=1}^k \rv_m (\sj F, F)
\end{CD}
$$
where $\Delta$ is the $k$-fold diagonal map.
\item $d_0 = \sum_{i=1}^k d_i : \rv_m(\sj F, F) \to  \rv(\sj^k F, F)$.
\end{enumerate}
\end{corollary}

\begin{proof}   By part (3) of   Proposition \ref{relations}, $e_i \circ d_0 = id : \rv_m(\sj F, F) \to \rv_m(\sj F, F).$   Thus $ \oplus_{i=1}^k e_i  \circ d_0$ is the diagonal map,
$\rv_m(\sj F, F) \to  \bigoplus_{i=1}^k \rv_m (\sj F, F)$.  This proves part (1).  For part 
(2), notice that by the above theorem, $\oplus_{i=1}^k e_i $ is an isomorphism, inverse to $\sum_{i=1}^k d_i$.   So by the commutativity of the diagram in part (1), $$d_0 = \left(\sum_{i=1}^k d_i\right)\circ \Delta,$$ 
which is to say, $d_0 (x) = \sum_{i=1}^k(d_i(x)).$
\end{proof}

\med
We next   consider an embedding,
$$
\iota : \sj^k F \hk \si^{k-1}\sj F
$$
defined as follows (see figures \ref{fig:figeight} and \ref{fig:fignine} below).   Let $P_k$ be the surface of genus zero with $k+1$ boundary components, as in figure \ref{fig:figeight}.  
\begin{figure}[ht]
  \centering
  \includegraphics[height=4cm]{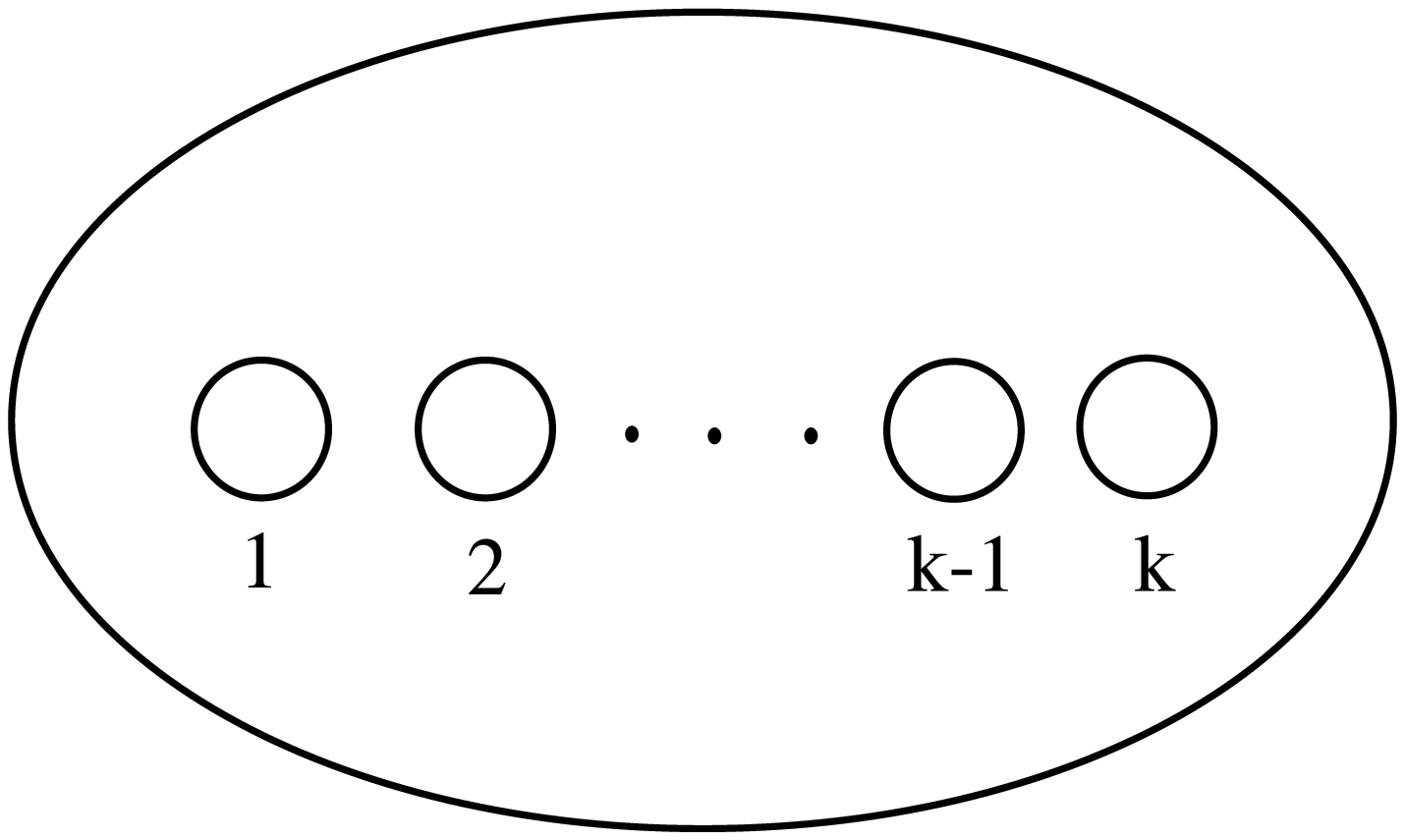}
 \caption{The surface $P_k$ } 
  \label{fig:figeight}
\end{figure}

$P_k$ is glued onto $\sj^k F$ by identifying the $k$- interior boundary circles of $P_k$ with those in $\sj^k F$ that have been created by the operation $\sj^k$.  See figure \ref{fig:fignine}.

\begin{figure}[ht]
  \centering
  \includegraphics[height=6cm]{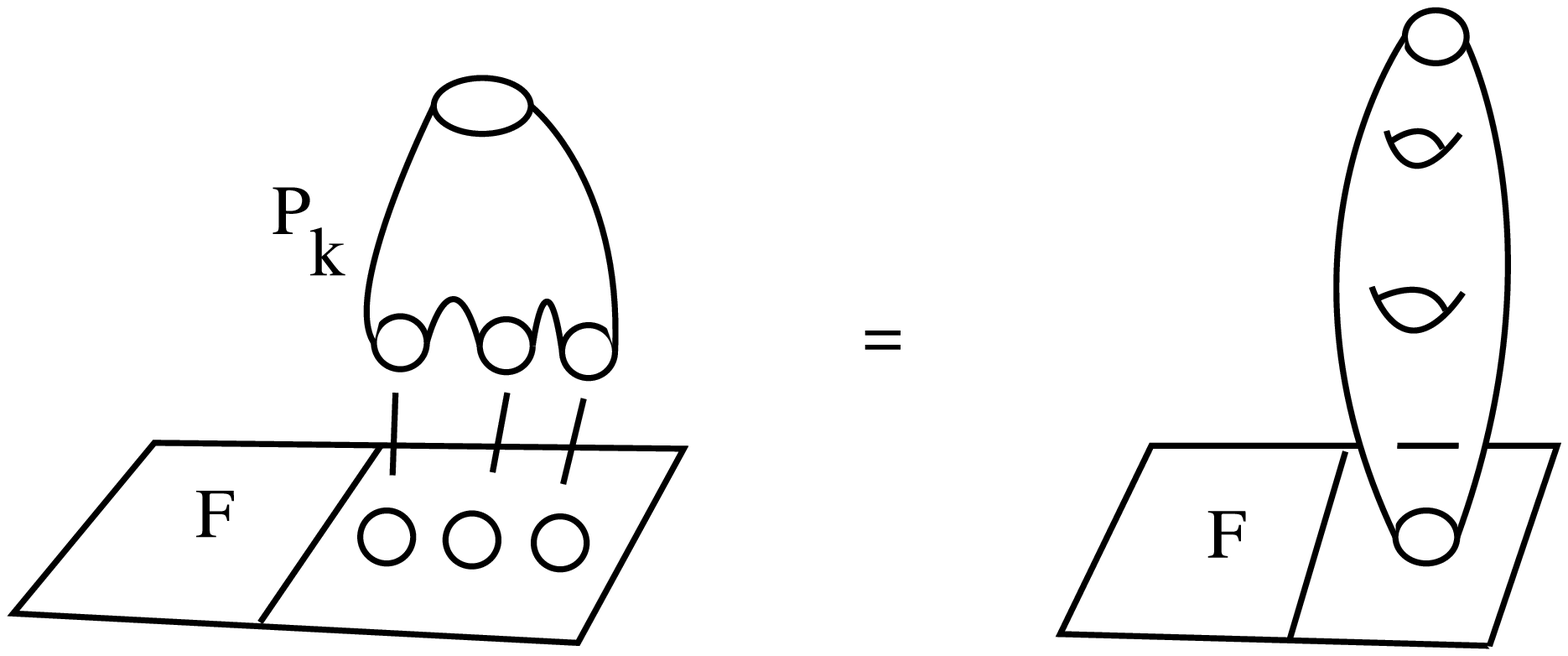}
 \caption{The embedding $\iota :  \sj^k F \hk \si^{k-1}\sj F$ } 
  \label{fig:fignine}
\end{figure}

Viewed as a morphism in the surface category $\cc$, $\iota$ induces a homomorphism,
$$
\iota_*: \rv_m(\sj^kF, F) \to \rv_m(\si^{k-1}\sj F, F).
$$
We now prove the following relations concerning compositions, $\rv_m(\sj F, F) \to \rv_m(\sj^kF, F) \to \rv_m(\si^{k-1}\sj F, F).$

 \begin{lemma}\label{morerelations}
\begin{align}
\iota_* \circ d_0 &= \si^{k-1} : \rv_m(\sj F, F) \to \rv_m(\si^{k-1}\sj F, F) \notag \\
\iota_* \circ d_j &= \iota_*\circ d_q \quad \text{for} \quad  j, q  = 1, \cdots, k \notag
\end{align}
\end{lemma}

\begin{proof} The first of these statements is immediate by the definitions of $d_0$ and $\iota$.  We therefore concentrate on the proof of the second statement. 

For each pair, $j, q = 1, \cdots , k$,  one can find  an element of the mapping class group,
$g_{j,q} \in \Gamma (\si^{k-1}\sj F)$ represented by a diffeomorphism  that is fixed on $\sj F$,   such that the induced embedding 
$$
g_{j,q} \circ \iota \circ d_j \quad \text{is isotopic to} \quad \iota \circ d_q : \sj F \hk \si^{k-1}\sj F.
$$
For example, $g_{j,q}$ can be taken to be the half Dehn twist around the curve $C_{j,q}$ depicted in figure \ref{fig:figten} below. 
 
 \begin{figure}[ht]
  \centering
  \includegraphics[height=6cm]{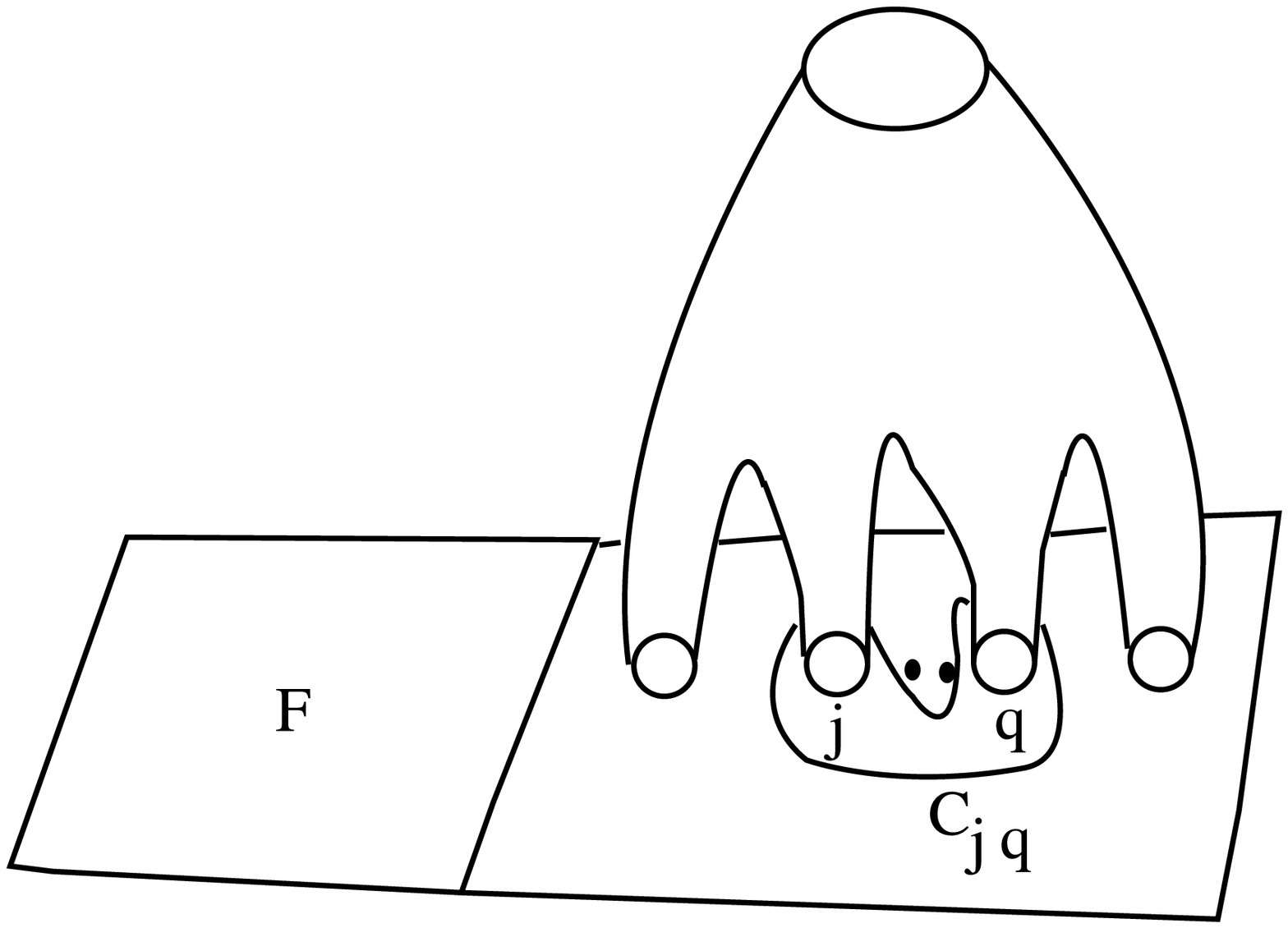}
 \caption{The  diffeomorphism $g_{j,q}$ given by the half Dehn twist around $C_{j,q}$} 
  \label{fig:figten}
\end{figure}

Now recall that given any embedding, $\phi : F_1 \hk F$, then there is an induced homomorphism of mapping class groups,
$$
\phi_* : \G (F_1) \to \G(F_2)
$$
defined by extending a diffeomorphism of $F_1$ to be the identity outside of the image of $\phi (F_1)$ in $F_2$.   
So in particular, for any $h \in \G(F_1)$, represented by a diffeomorphism,  there is a   diagram of embeddings that commutes up to isotopy,
$$
\begin{CD}
F_1  @>\phi >\hk >  F_2 \\
@Vh VV   @VV\phi_*(h) V \\
F_1  @>\phi >\hk >  F_2
\end{CD}
$$

This means that for any element of the mapping class group, $h \in \Gamma (\sj F)$, the following diagram of embeddings commutes up to isotopy.
$$
 \begin{CD}
 \sj F   @>\iota \circ d_j >\hk >  \si^{k-1}\sj F   @>g_{j,q} >\cong >   \si^{k-1}\sj F \\
 @Vh V \cong V    @V(\iota \circ d_j)_*(h) V\cong V   @V\cong V g_{j,q}(\iota \circ d_j)_*(h)g_{j,q}^{-1} V \\
  \sj F   @>\hk >\iota \circ d_j  >  \si^{k-1}\sj F   @>\cong>g_{j,q} >    \si^{k-1}\sj F
 \end{CD}
$$
Now the compositions in the  horizontal rows in this  diagram are isotopic to the  embedding $\iota \circ d_q$.
This means that the diffeomorphism  $g_{j,q}(\iota \circ d_j)_*(h)g_{j,q}^{-1}$  represents the class $(\iota \circ d_q)_*(h).$  In other words, the following diagram of homomorphisms of mapping class groups commutes:
$$
\begin{CD}
\G(\si^{k-1}\sj F)    @>\text{conjugation by} \,  g_{j,q} >> \G(\si^{k-1}\sj F) \\
@A(\iota \circ d_j)_* AA    @AA (\iota \circ d_q)_* A \\
\G(\sj F)   @>>= >   \G(\sj F)
\end{CD}
$$

But in homology of groups, conjugation by a group element acts as the identity.  So we have
$$
(\iota \circ d_j)_* = (\iota \circ d_q)_* : \rv_m (\sj F, F) \to \rv_m(\si^{k-1}\sj F, F) 
$$ as claimed. 
\end{proof}

\med
We now prove Proposition \ref{complete}. 
\begin{proof}
 We will show that for every $k > 0$,  there is a homomorphism,
$$
\psi_k : \rv_m(\sj F, F) \to \rv_m(\sj F, F)
$$
so that for any $\alpha \in \rv_m(\sj F, F)$,  $\alpha = k\cdot \psi_k (\alpha).$Since $\rv_m(\sj F, F)$ is a finitely generated abelian group, this  will imply that $\rv_m(\sj F, F) = 0$. 

Consider the homomorphism 
$$
\si^{k-1} : \rv_m(\sj F, F) \to \rv_m(\si^{k-1}\sj F, F).
$$
By Ivanov's theorem (\cite{ivanov}, Corollary 4.9),    $\si^{k-1}$ is an isomorphism for $m$ in the range $g \geq 2m+d+2$.    This allows us to compute:

\begin{align}
\alpha &= (\si^{k-1})^{-1} \circ (\si^{k-1})(\alpha) \notag \\
&= (\si^{k-1})^{-1}(\iota_*\circ d_0)(\alpha) \quad \text{by Theorem \ref{morerelations}} \notag \\
&=(\si^{k-1})^{-1}(\iota_*\circ \sum_{i=1}^k d_i(\alpha)) \quad \text{by part (2) of Corollary \ref{moredees}} \notag \\
&= (\si^{k-1})^{-1}(k\cdot \iota_*(d_1(\alpha)) \quad \text{by Theorem \ref{morerelations}} \notag \\
&=k\cdot ((\si^{k-1})^{-1}(\iota_*d_1(\alpha)). \notag
\end{align}

So we define 
$$
\psi_k = (\si^{k-1})^{-1}\circ \iota_*\circ d_1 : \rv_m(\sj F, F) \to \rv_m(\sj^k F, F) \to \rv_m(\si^{k-1}\sj F, F) \xr{\cong} \rv_m(\sj F, F).
$$
This completes the proof of the Proposition \ref{complete}, and hence completes the inductive step in the proof of Theorem \ref{mpg2} for the operation $\sj$.   \end{proof}

   \subsection{Closing the last hole }
   
   Our goal in this section is to prove Theorem \ref{mpg2} for the operation $\sk$.  This will complete the proof of Theorem  \ref{mpg2}, and therefore of Theorem \ref{mpg}.  We continue to operate under Inductive Assumption \ref{assume1}.  
   
    Now in the previous two sections, we proved that for any surface $F$ with at least one boundary component, $\rv_m(\sj F, F) = 0$ for $2m \leq g(F) - d+2$, where $g(F)$ is the genus. From Lemma \ref{unique} we have the following.     
   
   \begin{corollary} \label{twoboundary}  Let $F$ be a surface with at least two parameterized boundary components, and let $V$ be a coefficient system of overall degree $\leq d$.    Then
   $$\rv_m(\sk F, F)  = 0$$
   for $2m \leq g(F) - d+2$.  
   \end{corollary}

   Thus to complete the proof of Theorem \ref{mpg2}, we can restrict our attention to the case when $F$
   is  a surface of genus $g$ with one parameterized boundary component.  That is, $F$ is an object in the category $\cc_{g,1}$.  Let $\hf = \sk F$ be the associated closed surface.  The goal of this section is to prove the following theorem,  which generalizes  Theorem 1.9 of \cite{ivanov} to non-trivial coefficient systems. 
   
   \begin{theorem}\label{hole}  Let $V$ be a coefficient system of overall degree $d$.      Then
   $$
 \rv_m(\hf, F)  = 0$$
   for $2m \leq g(F) - d+1$.  
   \end{theorem}
   
   \med

    \begin{proof}        
   We will prove Theorem \ref{hole} by using Corollary \ref{twoboundary} and a set of spectral sequences.   An important tool that we will use is the well known Zeeman comparison
   theorem for spectral sequences \cite{zeeman}.   
      
   \med
   Consider two first quadrant spectral sequences, $E = \{E^r_{p,q}, d^r\}$ and $\bar E = 
   \{\bar E^r_{p,q}, \bar d^r\}$ converging to graded groups $A= \{A_n\}$ and $\bar A = \{\bar A_n\}$ respectively.  Let $f : E \to \bar E$ be a morphism of spectral sequences converging to a morphism $g : A \to \bar A$ of graded groups.

   We shall make use of the comparison theorem for the map $f$, and state the form we need.  We call a map $g : A \to \bar A$ of graded groups $N$-connected if $g_m : A_m \to \bar A_m$ is an isomorphism for $m < N$,  and an epimorphism for $m = N$. 
   
   Consider the following condition of the map $f : E \to \bar E$, often ensured by the universal coefficient theorem and  ``fiber information".  
   
 \bfl
 \bf Condition (*).  \sl 
  If $f_{p,0} : E^2_{p,0} \to \bar E^2_{p,0}$ is an isomorphism for $p \leq P$, and if $f_{0,q} : E^2_{0,q} \to \bar E^2_{0,q}$ is an isomorphism for $q \leq Q,$ then $f_{p,q} : E^2_{p,q} \to \bar E^2_{p,q}$ is an isomorphism for all $p \leq P$ and $q \leq Q$.
\rm

\efl
   
   \med

    \begin{theorem}\label{comparison}   \begin{enumerate}
    \item Suppose that $f^2_{p,q} : E^2_{p,q} \to \bar  E^2_{p,q}$ is an isomorphism for $p+q <N$ and an epimorphism for $p+q = N$.  Then $g : A \to \bar A $ is $N$-connected.
    \item If $g : A \to \bar A$ is $N$-connected, and condition $(*)$ is satisfied, then $f_{*,0} : E^2_{*,0} \to \bar E^2_{*,0}$ is $N$-connected.
    \end{enumerate}
       \end{theorem}
       
       \med
       In our application of Theorem \ref{comparison}, the spectral sequences come from filtered chain complexes, and  $f$ from a map of filtered chain complexes.  In this situation there is a relative spectral sequence 
       $ \{E^r_{p,q}(g)\}$ converging to $H_{p+q}(g)$, and an exact sequence
       $$
       \cdots \to E^2_{p,q} \xr{f} \bar E^2_{p,q} \to E^2_{p,q}(g) \to E^2_{p-1,q} \xr{f} \cdots
       $$
       that makes Theorem \ref{comparison} obvious.  This is J. Moore's original argument in \cite{cartan}

     We will now describe the spectral sequences we will be studying in order to prove   Theorem \ref{hole}.     
   Recall from   \cite{harer} and \cite{ivanov} the curve complex $C_0 (F)$.  A vertex in this complex is an isotopy class of a nontrivial\footnote{Nontrivial means that the circle $L$ cannot be deformed to a point or to a boundary circle.}  embedded circle $L$ so that the complement $F - L$ is connected.   A $p$-simplex $L^p = \{L_0, \cdots , L_p\}$ is a set of disjoint embedded circles subject to the condition that $F-L^p$, by which we mean
   $F - \bigcup_{i=0}^p L_i$, is connected.  Notice that if $F$ has genus $g$ and $r$-boundary components, $F_{g,r} - L^p \cong F_{g-p-1, 2p+2+r}$ has genus $g-p-1$ and $2p+2+r$ boundary components.    This implies that the complex $C_0(F_{g,r})$ has dimension $g-1$.   In fact, by Harer's Theorem  1.1 of \cite{harer},
   $$
   C_0(F_{g,r}) \simeq \bigvee S^{g-1}.
   $$
   
   We now consider the case $F = F_{g,1}$ so that $\hf $ is a closed surface of genus $g$. The embedding $j : F \hk \hf$ induces a map of curve complexes,
   $$
   j : C_0 (F) \la C_0(\hf)
   $$
   that is equivariant with respect to $\sk : \G(F) \to \G(\hf)$.
   
   \med
   We will study the relative spectral sequence (\ref{relspec}) for the map $j$ with respect to the group
   actions of $\G (F)$ and $\G(\hf)$.  
   
   \med
   Notice that for any surface $R$, the mapping class group $\G (R)$ acts transitively on the $p-1$-simplices of $C_0(R)$.  The isotropy group of a $(p-1)$-simplex $L$, $\G(R)_L$,  can permute the vertices of $L$.   This induces a surjective homomorphism to the symmetric group,
   $$
   p : \G(R)_L \to \Sigma_{p}.
   $$
   The relation between $\G(R)_L$ and $\G(R-L)$ is expressed in the following group extensions (\cite{ivanov}, section 1.8), the second of which is central:
   
   \begin{align}\label{onea}
   1 \to \tilde \G(R)_L \to &\G(R)_L \xr{p}  \Sigma_{p}\smallint \bz/2 \to 1 \\
   0 \to \bz^{p} \to &\G(R-L)  \to  \tilde \G(R)_L \to 1. \notag 
   \end{align}
   Here $\tilde \G(R)_L \subset \G(R)_L $ is the subgroup of mapping classes that fix each vertex of $L$ as well as their orientations. (See \cite{ivanov}, p. 159 for notation.)
   
   \med
   
   We continue to use the notation introduced in the description of the relative spectral sequence (\ref{relspec}) above.   
   
   \med
   The relative spectral sequence for the pair $(\hf, F)$  has $E^1$ -term given by
   \begin{equation}\label{two}
   E^1_{p,q}(\hf, F) = H_q(\G (\hf)_{j(L)}, \, \G(F)_{L};  V(\hf)_{j(L)}, V(F)_{L}).
   \end{equation}
   Here $V(F)_L = V(F) \otimes \bz_L$, where $\bz_L$ is the sign representation, $\G(F) \to \Sigma_p \to \{\pm 1\}$. 
      Notice that by the connectivity of the curve complex,  this spectral sequence converges to zero
in total degrees $\leq g$.     Notice furthermore that  on the edge we have,
$$
E^1_{0,q} = H_q(\G(\hf), \G (F); \, V(\hf), V(F)) = \rv_q(\hf, F).
$$

More generally we will use the exact sequences  (\ref{onea}) to relate $E^1_{p,q}$ to $\rv_q(\hf-L; \,F-L)$. Indeed the Hochschild-Serre spectral sequences for these exact sequences takes the form,

\begin{align} 
E^2_{r,s} =&H_r(\Sigma_{p }\smallint \bz/2; \, H_s(\tilde \G(F)_{L}; V(F)_{L}) \quad \text{converging to} \quad H_{r+s}(\G (F)_{L}; \, V(F)_{L}) \label{four} \\
E^2_{r,s} =&H_r(\tilde \G(F)_L; \, V(F)) \otimes H_s(\bz^{p}) \quad \text{converging to} \quad H_{r+s}( \G(F-L); \, V(F)) \label{five}
\end{align}

\med
Each of these spectral sequences maps to the corresponding spectral sequence for $\hf$, replacing $F$.  Our intent is to use the   comparison theorem for spectral sequences for these maps. 

\med
We first look at this map on the  second spectral sequence (\ref{five}).  On the $E^2$-level this is a map

$$
j_*: H_r(\tilde \G(F)_L, V(F)) \otimes H_s(\bz^p) \la H_r(\tilde \G(\hat F)_{j(L)}, V(\hf)) \otimes H_s(\bz^p)
$$
where $L$ represents an arbitrary $(p-1)$-simplex of $C_0(F)$.  We shall use Theorem \ref{comparison} part (2) to study its connectivity.

We write $F-L \, = \, \fgp$  as it is a surface of genus $g-p$ with $2p+1$-boundary components.  Moreover, $$\hf - j(L) = \sk (F - L). $$

The spectral sequences converge to $H_*(\G (F-L); V(F))$ and $H_*(\G(\hat F-j(L)); V(\hat F))$, respectively.  The map
\begin{equation}\label{six}
j_* : H_r(\G(F-L); V(F)) \to H_r(\G(\hat F-j(L)); V(\hat F))
\end{equation}
  fits into a long exact sequence with relative terms $$Rel_r^{\Sigma_{p, -2p}V}(\sk \fgp, \, \fgp). $$    

 By  Corollary \ref{twoboundary} and   Induction Assumption \ref{assume1} on the degree of the coefficient system,
$$
 Rel_r^{\Sigma_{p, -2p}V}(\sk \fgp, \, \fgp) \cong  Rel_r^{V}(\sk \fgp, \, \fgp)
 $$ 
 for $2r \leq g-p-d$.   But this group is zero by Corollary \ref{twoboundary}.     Thus the map in $(\ref{six})$ is $(g-p-d)/2$-connected.   
  
 By the spectral sequence comparison theorem \ref{comparison} part (2), it follows that
 
 \begin{equation}\label{seven}
 j_* : H_r(\tilde \G(F)_{L}; \, V(F)) \to H_r(\tilde \G(\hf)_{j(L)}; \, V(\hf))
 \end{equation}
 is  also $(g-p-d)/2$-connected.   
 
 Inputting this into the $E^2$-term of the first spectral sequence of (\ref{four}) and  using  comparison theorem \ref{comparison} part (1)  shows that 
 $$
 j_* : H_{*}(\G(F)_{L} ; \, V(F)_{L}) \to H_{*}(\G(\hf)_{j(L)}; \, V(\hf)_{j(L)})
 $$
 is $(g-p-d)/2$-connected. 
  
 \med
 Now we compare with (\ref{two}) and see that
 $$
 E^1_{p,q}(\hf, F) = 0 \quad \text{if $p \geq 1$ and $p+2q \leq g- d.  $}
 $$
Since $p+2q \leq 2(p+q)-1$ when $p \geq 1$, it follows that $E^1_{p,q}(\hf, F) = 0$ in total degree $2(p+q)\leq g - d+1$ for $p \geq 1$.  However we know this spectral sequence converges to zero in this range.  This implies that $E^1_{0,m}(\hat F, F) = 0$ for $2m \leq g- d-1$.  But by (\ref{two}) this is the relative homology group $\rv_m(\hf , F)$.   This completes the proof of Theorem \ref{hole}.  \end{proof}   

\med
We end by recalling that Theorem \ref{hole} completes the inductive step in the proof of Theorem \ref{mpg2} for the operation $\sk$, which was our last step in our proof of Theorem \ref{mpg2} and therefore of Theorem \ref{mpg}   of the introduction.  Indeed, notice that we have proved the following strengthening of Theorem \ref{mpg}.

Let $\phi : F_1 \hk F_2$ be an embedding of the sort used in defining a   morphism in the surface category $\cc$.  It induces
a homomorphism of mapping class groups $\phi : \G (F_1) \to \G(F_2)$ by extending a diffeomorphism  of $F_1$ that fixes its boundary, to a diffeomorphism of $F_2$ by letting it act as the identity on the complement $F_2 - \phi (F_1)$.    $\phi$ then induces
a homomorphism in homology, $\phi_*: H_q(\G (F_1); V(F_1)) \to H_q(\G (F_2); V(F_2))$. 

\med
\begin{theorem}\label{strong}
Let $\phi : F_1 \to F_2$ be any embedding that defines a morphism in the surface category $\cc$, where the genera of these surfaces are $g_1$ and $g_2$, respectively.  Then if $V$ is a coefficient system of degree $d$, the induced homomorphism
$$
\phi_* : H_q(\G (F_1); V(F_1)) \to H_q(\G (F_2); V(F_2))
$$
is an isomorphism if $2q+d+2 < g_1$, and is an epimorphism if $2q+d+2= g_1$.
\end{theorem}

\med
\begin{proof}.  First notice that for    such an embedding $\phi$ to exist,    we must have $g_2 \geq g_1$.  The theorem now follows because
any morphism  is isotopic to  a composition of the embeddings $e_{i,j} : F \hk \Sigma_{i,j}F$ for $(i,j)$ of the form $(1,0)$, $(0,1)$, or $(0, -1)$,  as well as diffeomorphisms.  The work in this section implies this result about each of these morphisms, and hence about any composition of these morphisms. \end{proof}

   \section{Stability of the space of surfaces}
   
   Our goal in this section is to prove the first part of Theorem  \ref{unstable}, as stated in the introduction.
   This is a stability theorem for the surface spaces, $\sgn(X; \gamma)$.  We will prove this theorem in two steps.  In order to describe these steps it is helpful to introduce
   the following terminology.
   
   \begin{definition}\label{algstab}  We say that a space $X$ is algebraically stable if for every surface $F_{g,n}$  having genus $g$ and $n$  boundary components,  the homology groups,
   $$
   H_p(\G_{g,n}; H_q(Map_\p (F_{g,n}, X)))
   $$
   are independent of $g$, $n$,    for $2 p+ q  < g $.  By this we mean that if   $\phi : F_1 \to F_2$ is any embedding that defines a morphism in the surface category $\cc$, the induced map
  $$
  \phi_* : H_p(\G (F_1),  H_q(Map_\p (F_1, X))) \to H_p(\G (F_2), H_q(Map_\p (F_2, X)))
  $$
  is an isomorphism so long as the genus $g(F_1) > 2p+q+2 $, and is an epimorphism if $g(F_1) =  2p+q+2.$   Here $Map_\p(F, X)$
  refers to those maps that send the boundary $\p F$ to a fixed basepoint in $X$.
 \end{definition}
 
 \begin{lemma}\label{emstable} If $X$ is an Eilenberg-MacLane space, $X = K(A,m)$ with $m \geq 2$,
 then $X$ is algebraically stable.
 \end{lemma}
 
 \begin{proof} This follows from  Example 2 given after the definition of the degree of a  coefficient system (Definition \ref{degree}) in section 1.1, and from Theorem \ref{strong} above.  
  \end{proof}  
 
 \med
 The two steps that we will use to prove Theorem  \ref{unstable} are the following.

 \med
 \begin{theorem}\label{xstab}  Every simply connected space is algebraically stable.
 \end{theorem}

 \med
 \begin{theorem}\label{topstab}  If $X$ is algebraically stable, then the homology
 of the surface space, $H_p(\sgn (X; \gamma))$ is independent of $g, n$ and $\gamma$ if $2p+4 \leq g$.  
  \end{theorem}

 \med 
\sl Proof of Theorem \ref{xstab}.   \rm

    To prove this theorem we use a classical tool of obstruction theory, the Postnikov tower of a simply connected space $X$.  This is a sequence of maps 
  $$
 \begin{CD}
 \cdots X_m @>>>  X_{m-1} @>>> \cdots @>>> X_3 \to X_2= K(\pi_2(X), 2) \\
 && @Vk_{m-1}VV  && @Vk_2VV \\
 && K(G_m , m+1)   &&&&  K(G_3 , 4)
 \end{CD}
 $$ where the $G_i$'s are the homotopy groups, $G_i = \pi_{i}(X)$,   and the spaces $K(G_i, i+1)$ are Eilenberg-MacLane spaces.  This tower satisfies   the following properties.  
 \begin{enumerate}
 \item  Each map $  X_{m-1} \xr{k_{m-1}} K(G_m, m+1)$ is a fibration  with fiber $X_{m}$.  
 \item  The tower comes equipped with maps $f_i : X \to X_i$ which are $(i+1)$-connected.   
 \end{enumerate}
 
 Let $F_{g,n}$ be a fixed surface in $\cc_{g,n}$.    
  The above Postnikov tower gives rise to an induced tower,
 
  $$
 \begin{CD}
 \cdots \to  Map_\p(F_{g,n}, X_m)  @>>>  Map_\p(F_{g,n},X_{m-1})  @>>> \cdots @>>>  Map_\p(F_{g,n},  K(\pi_2(X),2) ) \\
 && @Vk_{m-1}VV  && @Vk_2VV \\
 && Map_\p(F_{g,n},   K(G_m, m+1) )   &&&& Map((S,  K(G_3 , 4))
 \end{CD}
 $$ where each  map $k_{m-1} : Map_\p(F_{g,n},X_{m-1}) \to Map_\p(F_{g,n},   K(G_m, m+1) )$  is a  fibration with fiber $ Map_\p(F_{g,n}, X_m)$.    This tower converges to $Map_\p (F_{g,n}, X)$.  

 By lemma \ref{emstable},   we know that $X_2 = K(\pi_2(X), 2)$ is algebraically stable.  We inductively assume $X_j$ is algebraically stable for $j \leq m-1$.  We now study $X_m$ by analyzing the homotopy fibration sequence,
 
 \begin{equation}\label{fib}
 Map_\p(F_{g,n}, K(G_m, m )) \to  Map_\p(F_{g,n}, X_m) \to Map_\p(F_{g,n},X_{m-1}).
\end{equation}

Notice that $H_p(\G (F_{g,n}); H_s (Map_\p (F_{g,n},X_{m-1}))  \otimes H_t(Map_\p(F_{g,n},   K(G_m, m)))$ is independent of $g$ and $n$ so long as $2p+s+t < g-2.$  The coefficients are the $E^2_{s,t}$- term of the Serre spectral sequence for this fibration.  Notice furthermore that this spectral sequence   of coefficient systems.  That is,  for each $r$,  $ F_{g,n} \to E^r_{p,q}(F_{g,n})$
  is a coefficient system, and the differentials are natural transformations.   We therefore know that $H_p(\G (F), E^2_{s,t}(F))$ is independent of the surface $F$, so long as the genus $g(F)$ satisfies, $g(F)-2 > 2p+s+t$.  This means that if $\phi : F_1 \to F_2$ is an embedding that defines a morphism in the surface category $\cc$, the induced map
  $$
  \phi_* : H_p(\G (F_1), E^2_{s,t}(F_1)) \to H_p(\G (F_2), E^2_{s,t}(F_2))
  $$
  is an isomorphism so long as $2p+s+t < g(F_1)-2, $  and an epimorphism if $2p+s+t = g(F_1)-2$.   Inductively assume the same statement is true for $E^r_{s,t}$ replacing $E^2_{s,t}$ in this mapping.  The fact that the differential $d_r : E^r_{s,t} \to E^r_{s-r, t+r-1}$ is a natural transformation, means the following diagram commutes:   
  $$
  \begin{CD}
   H_p(\G (F_1), E^r_{s,t}(F_1))     @>\phi_*>\cong > H_p(\G (F_2), E^r_{s,t}(F_2)) \\
   @Vd_r VV   @VV d_r V \\
   H_p(\G (F_1), E^r_{s-r,t+r-1}(F_1))     @>\phi_*>\cong > H_p(\G (F_2), E^r_{s-r,t+r-1}(F_2))
   \end{CD}
   $$
   Passing to homology we may conclude that for $g(F_1)-2 >  2p+s+t $,  the map
  $$
  \phi_* :   H_p(\G (F_1), E^r_{s,t}(F_1))     \to H_p(\G (F_2), E^r_{s,t}(F_2))
  $$
  is an isomorphism, and is an epimorphism if $g(F_1)-2 = 2p+s+t$.    We may therefore conclude that the same statement holds when the coefficients are the $E^\infty$-term.  That is,
  $$
 \phi_* :   H_p(\G (F_1), H_{s+t}(Map_\p(F_1, X_m))   \to H_p(\G (F_2),(Map_\p(F_2, X_m)) 
 $$
 is an isomorphism if $g(F_1)-2 >  2p+s+t$, and an epimorphism if $g(F_1)-2 = 2p+s+t$.  In other words, $X_m$ is algebraically stable.  By induction we may conclude that $X$ is algebraically stable.
\hfill \openbox

 \bg
   
   We now complete the proof of the first part of  Theorem \ref{unstable} by proving Theorem \ref{topstab}.  
   
   \begin{proof}  Let $X$ be an algebraically stable space.  We need to prove that
   $H_p(\sgn(X; \gamma))$ does not depend on $g$ or $n$ if $2p+2 \leq g-2$.  
   
     Let $\Sigma$ denote one of the operations $\si$, $\sj$, or $\sk$.  It suffices to show that $\Sigma$ induces a homology isomorphism, $$\Sigma_* : H_p(\sgn(X; \gamma)) \xr{\cong} H_p(\msc S_{g+i, n+j}(X; \gamma))$$ for $2p+2 \leq g-2$.  Here $i$ and $j$ depend in the obvious way on which operation ($\si$, $\sj$, or $\sk$) $\Sigma$ represents. 
   
    Recall that the homotopy type of $\sgn (X; \gamma)$ does not depend on the choice of $\gamma$.
   So we may assume that $\gamma : \coprod_n S^1 \to X$ is the constant map.  
   Consider the fibration described in the introduction
$$
Map_\p(F_{g,n}, X)  \to\sgn(X; \gamma) \to \sgn (point) =  BDiff^+(F_{g,n}, \p).
$$

    The map $\Sigma$ 
  induces a map  from this fibration
 to the corresponding fibration where $F_{g,n}$ is replaced by $\Sigma (F_{g,n})$.    They therefore induce maps of the corresponding Serre spectral sequences
 $$
\Sigma_* : E^r_{p,q}(F_{g,n}) \to E^r_{p,q}(\Sigma (F_{g,n}).
 $$

  On the $E^2$ level these are homomorphisms
\begin{align}
\Sigma_* : &H_p(BDiff(F_{g,n}, \p  ); \{H_q (Map_\p(F_{g,n},    X  )\}) \notag \\
   & \la H_p(BDiff(\Sigma(F_{g,n}), \p  ); \{H_q (Map_\p(\Sigma (F_{g,n}),    X  )\}).\notag
  \end{align} where $\{\}$ indicates that $\pi_1(BDiff(F_{g,n},\p))=\G(F)$ acts on $H_q(Map_\p(F_{g,n},    X  ))$.
  
Now since these diffeomorphism groups are homotopy discrete, that is, each path component of these groups is contractible, and since the corresponding groups of path components are the   mapping class groups, this homomorphism is given
by
 
$$
\Sigma_* : H_p(\Gamma (F_{g,n}),  \{H_q (Map_\p(F_{g,n},    X  )\}) \to H_p(\Gamma (\Sigma (F_{g,n}), \{H_q (Map_\p(\Sigma (F_{g,n}),    X  )\}).
$$

Now under the assumption that $X$ is algebraically stable, this map is an isomorphism for $2p+q < g-2$, and an epimorphism for $2p+q=g-2$.       That is, the homomorphisms

  $$
\Sigma_* : E^2_{p,q}(F_{g,n})) \to E^2_{p,q}(\Sigma (F_{g,n}))
 $$
are isomorphisms   for $2p+q  < g-2$, and    epimorphisms for $2p+q=g-2$.
   By the Zeeman comparison theorem \ref{comparison},  this implies that on the $E^\infty$ level, the maps
   
     $$
\Sigma_* : E^\infty_{p,q}(F_{g,n}) \to E^\infty_{p,q}(\Sigma (F_{g,n}))
 $$
are isomorphisms for $2p+q < g-2$.     The theorem
now follows by the convergence of the Serre spectral sequence.
 \end{proof}

   \section{The stable topology of the space of surfaces, $\sgn (X; \gamma)$}
   Our goal in  this section is to prove Theorem \ref{one} as stated in the introduction.
   Our method is to use the results of   \cite{GMTW} on cobordism categories and to adapt the methods of McDuff-Segal \cite{mcduffsegal} and Tillmann \cite{tillmann} on group completions of categories.  Alternatively, one could use the argument given in \cite{madsenweiss}, section 7. 
   
   \subsection{The cobordism category of oriented surfaces mapping to $X$}
   
   The topology of cobordism categories was described in great generality in \cite{GMTW}.  We describe their result as it pertains to our situation.  
   
   \med
   \begin{definition}\label{cat}
   Let $X$ be a simply connected, based space.  Define the category  $\cc_X$  of surfaces mapping to $X$, as follows.  
 
 The objects of $\cc_X$ are given by pairs $(C, \phi)$, where $C$ is a closed, oriented one-dimensional manifold, properly embedded in infinite dimensional Euclidean space, $C \subset \br^\infty$.   $\phi : C \to X$ is a continuous map.

 A morphism $(S, \psi)$ from $(C_1, \phi_1)$ to $(C_2, \phi_2)$ consists of an oriented surface $S$    properly embedded  $S \subset \br^\infty \times [a,b]$ for some interval $[a,b]$, together with a continuous map $\psi : S \to X$.  
 \end{definition}
 
 In this description, ``properly" embedded means  the following.
 The boundary $\p S$  lies in the boundary $\br^\infty \times \{a\} \sqcup \br^\infty \times \{b\}$.  The intersection of $S$ with these ``walls" are also assumed to be orthogonal in the sense that for sufficiently small $\eps >0$, the intersection of $S$ with
 $\br^\infty \times \{a+t\}$ and $\br^\infty \times \{b-t\}$ is constant for $0\leq t < \eps$.
 
 We write $\p_a S$ and $\p_bS$ for the intersection of $S$ with $\br^\infty \times \{a\}$ and $\br^\infty \times \{b\}$, respectively.  These are the ``incoming and outgoing" boundary
 components of $S$.   $(S, \psi)$ is a morphism from $(\p_a S, \psi_{|_{\p_aS}})$ to 
 $(\p_b S, \psi_{|_{\p_bS}})$
 
 Composition in this category is given by union of surfaces along common, parameterized boundary components.  In particular, assume  $(S_1, \psi_1)$ and $(S_2,\psi_2)$ are morphisms, with $S_1 \subset \br^\infty \times [a,b]$ and $S_2 \subset \br^\infty \times [c,d]$ with the target object of $(S_1, \psi_1)$ equal to the source object of $(S_2, \psi_2)$.  Then the composition is the glued cobordism, $(S_1 \#S_2, \psi_1 \#\psi_2)$, where $S_1 \#S_2  = S_1 \cup_{\p_b S_1 = \p_c S_2} S_2$, and $S_1 \#S_2 \subset \br^\infty \times [a, b+d-c]$.  $\psi_1 \#\psi_2 : S_1 \#S_2 \to X$ is equal to $\psi_1$ on $S_1$, and to $\psi_2$ on $S_2$.
 
 Finally we observe that $\cc_X$ is a topological category, where the objects and morphisms are topologized as described in the introduction.
 
 \med
 Let $|\cc_X|$ denote the geometric realization of the nerve of the category $\cc_X$. (This is sometimes called the classifying space of the category.)  The following was proved in \cite{madsenweiss}, but it is part of a more general theorem about cobordism categories proved in \cite{GMTW}.
 
 \begin{theorem}\label{cob1}
 There is a natural homotopy equivalence, $\alpha : \Omega |\cc_X| \xr{\simeq} \Omega^\infty (\bc \bp^\infty_{-1} \wedge X_+)$.
 \end{theorem}
 
 The right hand side is the zero space of the Thom spectrum $\bc \bp^\infty_{-1} \wedge X_+$ described in the introduction.  This is the infinite loop space appearing in
 the statement of Theorem \ref{one}.  Our approach to proving   Theorem \ref{cob1}  is to 
 use the stability result (Theorem \ref{unstable}) and a group completion argument to show
 how the stable surface space $\sinf (X; \gamma)$ is related to $\Omega |\cc_X|$. 
 First, however, for technical reasons, we need to replace $\cc_X$ by a slightly smaller cobordism category.  
 
 \begin{definition} Define the subcategory $\cred \subset \cc_X$ to have the same objects as $\cc_X$, but $(S, \psi)$ is a morphism in $\cred$ only if each connected component of $S$ has a nonempty outgoing boundary.
 \end{definition}
 
 The following was proved in \cite{GMTW}.
 
 \begin{theorem}\label{reduced} The inclusion $\cred \hk \cc_X$ induces a homotopy equivalence on geometric realizations,
 $$
 |\cred| \xr{\simeq} |\cc_X|.
 $$
 \end{theorem}
 
 Because of these two theorems, we have a homotopy equivalence, $\Omega |\cred| \xr{\simeq} \Omega^\infty (\bc \bp^\infty_{-1} \wedge X_+)$.   So to prove Theorem \ref{one}
 it suffices to prove the following.
 
 \begin{theorem}\label{groupcomplete} There is a map
 $$
 \beta : \bz \times \sinf (X; \gamma) \to \Omega |\cred|
 $$
 that induces an isomorphism in homology.
 \end{theorem}

   \subsection{A group completion argument and a proof of Theorem \ref{one}}
   
   We now proceed with a proof of Theorem \ref{groupcomplete}, which as observed above, implies Theorem \ref{one}.
   
   \begin{proof} Consider the following fixed object, $\sst$ of $\cred$.  $S^1$ is the unit circle in $\br^2 \subset \br^\infty$.  $* : S^1 \to x_0 \in X$ is the constant map at the basepoint. 
   
    For an object $(C, \gamma) $ of $\cred$, consider the space of morphisms, $Mor((C,\gamma), (\sst))$.  Suppose $C$ has $n-1$ components (i.e it is the union of $n-1$ circles embedded in $\br^\infty$).  Let $\psi : \sqcup_{n-1}S^1 \xr{\cong} C$ be a fixed parameterization, and let $\tilde \gamma :   \sqcup_{n-1}S^1  \to X$ be the composition,
    $\gamma \circ \psi$.  Finally let $\gamma^+ : \sqcup_n S^1 \to X$ be defined as follows.  Number the circles $0, \cdots, n-1$, and let $\gamma^+$  be  equal to the constant map at $x_0 \in X$ on the $0^{th}$ circle, and equal to $\tilde \gamma$ on    circles $1$ through $n-1$.      
    By definition of the morphisms in $\cred$, the following is immediate.
    
    \begin{lemma}\label{morph}  The morphism space  $Mor((C,\gamma), (\sst))$ is given by 
    $$
     Mor((C,\gamma), (\sst))  = \coprod_{g=0}^\infty \sgn(X; \gamma^+).
     $$
     \end{lemma}
    Now consider the morphism $\tst \in Mor (\sst, \sst)$, where $T$ is the surface of genus one described in the introduction, and $* : T \to x_0 \in X$ is the constant map.  
    For an object $(C,\gamma)$, define $Mor_{\infty}(C, \gamma)$ to be the homotopy colimit (or infinite mapping cylinder) under composing with the morphism, $\tst$,
    $$
    Mor_\infty (C, \gamma) = \hocolim \{Mor((C, \gamma), \sst)\xr{\circ \tst} Mor((C, \gamma), \sst) \xr{\circ \tst} \cdots \}
$$
An immediate corollary of the above lemma is the following.

\med
\begin{corollary}\label{morinf}
$$
Mor_\infty (C, \gamma) = \bz \times \sinf (X, \gamma^+).
$$
\end{corollary}

\med
Notice that $Mor_\infty$ is a contravariant functor,
\begin{align}
Mor_\infty : \cred &\to Spaces \notag \\
(C, \gamma) &\to Mor_\infty (C, \gamma). \notag
\end{align}
On the level of morphisms, if $(F, \phi)$, is a morphism from $(C_1, \gamma_1)$ to $(C_2, \gamma_2)$, the induced map 
$$
(F, \phi)^* : Mor_\infty (C_2, \gamma_2)  \to  Mor_\infty (C_1, \gamma_1)
$$
is given by precomposing with the morphism $(F, \phi)$.    Observe that Theorem \ref{unstable} implies the following.

\begin{lemma}\label{moriso}
 Every morphism in $\cred$, $$(F, \phi) : (C_1, \gamma_1) \to (C_2, \gamma_2)$$ induces a homology isomorphism
 $$
 (F, \phi)^* : Mor_\infty (C_2, \gamma_2)  \xr{\cong_{H_*}}  Mor_\infty (C_1, \gamma_1).
 $$
 \end{lemma}
 
 \med
 Now let $\cred \smallint Mor_\infty$ be the homotopy colimit of the functor $Mor_\infty$.  This is sometimes called the ``Grothendieck construction", and is modeled by the two sided bar construction, $B(*, \cred, Mor_\infty)$.  See \cite{mcduffsegal} and \cite{tillmann} for the details of this construction.  A consequence of Lemma \ref{moriso}, proved in \cite{mcduffsegal} is the following.
 
 \begin{proposition}\label{hfib}
 The natural projection map
 $$
 p : \cred \smallint Mor_\infty \la |\cred|
 $$
 is a homology fibration.  That is, the fibers of $p$ are homology equivalent to the homotopy fiber of $p$.
\end{proposition}   
We remark that this proposition, via Lemma \ref{moriso}, is the main place the theorems about the stability of the homology of mapping class group (Theorem \ref{mpg}) and the stability of the surface
space (Theorem \ref{unstable}) are used in the proof of Theorem \ref{one}. 

To complete   our proof of Theorem \ref{groupcomplete}, we need the following
result about the homotopy colimit space.

\begin{proposition}\label{contract}
The homotopy colimit space, $\cred \smallint Mor_\infty$ is contractible.
\end{proposition}

\begin{proof}  Consider the contravariant functor, $Mor_1 : \cred \to Spaces$ that assigns to an object, $(C, \gamma)$ the morphism space,
$$
Mor_1 (C, \gamma) = Mor((C, \gamma), \sst).
$$
Again on morphisms, $(F, \phi)$, $Mor_1$ acts by precomposition.   We observe that the homotopy colimit of this functor,
$ 
\cred \smallint Mor_1
$ is the geometric realization of the category whose objects are morphisms
in $\cred$ whose target is $\sst$, and where a morphism from $(F_1, \phi_1): (C_1, \gamma_1) \to \sst$ to $(F_2, \phi_2) : (C_2, \gamma_2) \to \sst$ is a morphism in $\cred$,  $(F, \phi) : (C_1, \gamma_1) \to (C_2, \gamma_2)$, so that
$$(F_2, \phi_2) \circ (F,\phi) = (F_1, \phi_1) : (C_1, \gamma_1) \to \sst.$$
But this category has a terminal object, $id : \sst \to \sst$, and hence its geometric
realization, $\cred \smallint Mor_1$ is contractible. 

Now composition with the morphism $\tst : \sst \to \sst$ defines a map,
$$
t : \cred \smallint Mor_1 \to  \cred \smallint Mor_1 ,
$$
and $\cred \smallint Mor_\infty$ is the homotopy colimit of the application of this map,
$$
\cred \smallint Mor_\infty =  \hocolim \{ \cred \smallint Mor_1 \xr{t}  \cred \smallint Mor_1 \xr{t} \cdots \}.
$$
Since it is a homotopy colimit of maps between contractible spaces, $\cred \smallint Mor_{\infty}$ is contractible.
\end{proof}
 We now complete the proof of Theorem   \ref{groupcomplete}.    By Propositions \ref{hfib} and \ref{contract}, we know that the map $p : \cred \smallint Mor_\infty \to |\cred|$
 has the property that if $(C, \gamma)$ is any object in $\cred$, which represents
 a vertex in $|\cred|$,  then there is a homology equivalence,
 $$
 p^{-1}(C, \gamma) \to \Omega |\cred|.
 $$
But by definition, $ p^{-1}(C, \gamma) = Mor_\infty (C, \gamma)$, which by Corollary \ref{morinf} is $\bz \times \sinf (X; \gamma^+).$     Thus we have a map
$$
\beta :  \bz \times \sinf (X; \gamma^+) \to \Omega |\cred|
$$
which is a homology equivalence.  But since, as we pointed out earlier, the homotopy
type of $\sinf (X; \alpha)$ does not depend on the boundary condition $\alpha$, Theorem 
\ref{groupcomplete} is proved.   As observed above, this completes the proof of Theorem \ref{one}.
   \end{proof}

 \end{document}